\documentclass[12pt]{article}
\usepackage{amssymb}
\usepackage{amsmath}

\input{tcilatex}
\begin{document}

\begin{center}
{\LARGE Self-similar vector measures of Markov-type operators}

\bigskip

Ion Chi\c{t}escu, Loredana Ioana, Radu Miculescu, Lucian Ni\c{t}\u{a}

\bigskip
\end{center}

\textbf{Abstract. }We consider iterated function systems (finite or
countable), together with linear and continuous operators on Hilbert spaces,
which enable us to construct Markov-type operators. Under suitable
conditions, these Markov-type operators have fixed points, which are
self-similar (invariant) vector measures, thus generalizing the classic
Hutchinson self-similar measures. Several models with concrete computations
are introduced.

\medskip



\bigskip

\textbf{Introduction}

\bigskip

In the present paper we introduce a generalization of self-similar measures
of the Markov operators generated by the Hutchinson construction which uses
iterated function systems with probabilities (see the seminal paper [12] and
[2], [9]). A description of this generalization follows. We replace the
positive measures (probabilities) in the classic model with vector measures
taking values in a Hilbert space $X$ and we replace the probabilities
forming the probability distribution with linear and continuous operators on 
$X$. We obtain operators in the space of variation bounded $X$-valued
measures, these operators having fixed points, given by the contraction
principle. The fixed points are self-similar (invariant) vector measures (we
call them also fractal measures). Of course, in order to use the contraction
principle, one must have complete metrics on the subsets of variation
bounded measures we are working with. These complete metrics are furnished\
by suitable norms, some of them being introduced in [6], using the integral
introduced in [5].

It is our duty to underline the strong influence of the paper [16]. Other
papers following similar lines are, e.g., [3], [10], [14] and [15] (more
closed to the ideas in the present paper being [3] and [15]). The recent
book [19] deals with the case of countable iterated function systems,
containing a large reference list. Other generalizations of the classic
Hutchinson construction can be found in [17] and [18]. We also mention the
study in [4], describing the influence of the measure $\mu $ upon the
support of the self-similar measure $\mu ^{\ast }$ (see the $H_{2}$-model).

The paper is divided as follows: introduction, three paragraphs and
references.

In the first paragraph we introduce the notations, notions and results used
throughout the paper. In order to make the paper self-contained, we briefly
recalled the essential contents of the papers [5] and [6].

The second paragraph forms the main part of the paper and is divided in
three subparagraphs. The first subparagraph introduces the theoretical
support and the underlying constructions (models). The second subparagraph
illustrates the theory with fixed point theorems in the space of vector
measures, accompanied by concrete examples. The third subparagraph is
dedicated to the particular case when all the contractions involved are
constant.

The third paragraph (divided in two subparagraphs) deals with the countable
case (the underlying generalized iterated function system has countably many
functions). Proofs are skipped or merely sketched, because they are similar
to those in the second paragraph. The first subparagraph introduces the
theory, while the second subparagraph introduces the underlying fixed point
theorems.

We mention that the essential contents of the present paper has been
presented by the first author at the 10th AIMS International Conference on
Dynamical Systems, Differential Equations and Applications, Madrid, 2014
(see Abstracts volume, pag 475).

We think that the reason for studying Markov-type operators on vector
measures and their fixed points is not only purely theoretical. Many
phenomena, e.g. behavior of fluids or of electric (magnetic) field are more
suited to a vector description.

It is our intention to continue the study of the present paper into two
directions: new and more variate applications (first direction) and
generalizations of the present constructions (second direction).

\bigskip

\textbf{1. Preliminary part}

\bigskip

Throughout this paper: $\mathbb{N}=\{0,1,2,...,n,...\}$, $\mathbb{N}^{\ast
}=\{1,2,...,n,...\}$, $\mathbb{R}_{+}=[0,\infty )$, $\overline{\mathbb{R}_{+}%
}=\mathbb{R}_{+}\cup \{\infty \}$ and $K$ will be the scalar field (real if $%
K=\mathbb{R}$, or complex if $K=\mathbb{C}$). All the sequences will be
indexed by $\mathbb{N}$ or $\mathbb{N}^{\ast }$ and all the vector spaces
(which are assumed to be non null) will be over $K$.

For any set $T$, $\mathcal{P}(T)$ is the set of all subsets of $T$. If $%
A\subset T$, $\varphi _{A}:T\rightarrow K$ is the characteristic (indicator)
function of $A$. If $T$ is a nonempty set, $X$ is a vector space, $\varphi
:T\rightarrow K$ and $f:T\rightarrow X$, we can consider the function $%
\varphi f:T\rightarrow X$ defined via $(\varphi f)(t)=\varphi (t)f(t)$, for
any $t\in T$ (many times, $f$ will be constant). The identity function $%
Id_{T}:T\rightarrow T$ acts via $Id_{T}(t)=t$, for any $t\in T$.

If $(E,\left\Vert .\right\Vert )$ and $(F,\left\vert \left\Vert .\right\Vert
\right\vert )$ are normed spaces, we consider the vector space $\mathcal{L}%
(E,F)=\{V:E\rightarrow F\mid V$ is linear and continuous$\}$ normed with the
operator norm $\left\Vert V\right\Vert _{o}=\sup \{\left\vert \left\Vert
V(x)\right\Vert \right\vert \mid x\in E,\left\Vert x\right\Vert \leq 1\}$
(which is even a Banach space if $F$ is a Banach space). In case $E=F$ we
write $\mathcal{L}(E)$ instead of $\mathcal{L}(E,E)$. If $F=K$, we write $%
E^{^{\prime }}$ instead of $\mathcal{L}(E,K)$ ($E^{^{\prime }}$ is the dual
of $E$). Considering the normed space $(E,\left\Vert .\right\Vert )$ (many
times we write only $E$), we have the weak$^{\ast \text{ }}$topology $\sigma
(E^{^{\prime }},E)$ on $E^{^{\prime }}$ (given by the family of seminorms $%
(\pi _{x})_{x\in E}$, where $\pi _{x}(x^{^{\prime }})=\left\vert x^{^{\prime
}}(x)\right\vert $, $x^{^{\prime }}\in E^{^{\prime }}$).

The scalar product of two elements $x,y$ in a Hilbert space $X$ will be
denoted by $<x,y>$. In case $X=K^{n}$, we have, for $%
x=(x_{1},x_{2},...,x_{n})$ and $y=(y_{1},y_{2},...,y_{n})$, the standard
scalar product $<x,y>=\underset{i=1}{\overset{n}{\sum }}x_{i}\overline{y_{i}}
$, generating the euclidean norm $\left\Vert x\right\Vert =(\overset{\infty }%
{\underset{n=1}{\sum }}\left\vert x_{n}^{2}\right\vert )^{\frac{1}{2}}$. For
a general Hilbert space $X$ with the scalar product $<.,.>$ and for $V\in 
\mathcal{L}(X)$, the adjoint of $V$ is $V^{\ast }\in \mathcal{L}(X)$ (hence $%
<V(x),y>=<x,V^{\ast }(y)>$ for any $x,y$ in $X$).

For any non empty set $T$ and any normed space $(X,\left\Vert .\right\Vert )$%
, we can consider the Banach space 
\begin{equation*}
B(T,X)=\{f:T\rightarrow X\mid f\text{ is bounded}\}
\end{equation*}%
equipped with the norm $f\mapsto \left\Vert f\right\Vert _{\infty }=\sup
\{\left\Vert f(t)\right\Vert \mid t\in T\}$ (the norm of uniform
convergence). We shall work in the particular situation when $(T,d)$ is a
compact metric space ($T$ having at least two elements). Then we have $%
C(T,X)\subset B(T,X)$, where $C(T,X)=\{f:T\rightarrow X\mid f$ is continuous$%
\}$ is a Banach space when equipped with the induced norm $\left\Vert
.\right\Vert _{\infty }$. Many times we write only $C(X)$ (resp. $B(X)$)
instead of $C(T,X)$ (resp. $B(T,X)$).

Let $(T,d)$ and $(X,\rho )$ be two metric spaces, $T$ having at least two
elements and let $f:T\rightarrow X$. The Lipschitz constant of $f$ is
defined by the formula%
\begin{equation*}
\left\Vert f\right\Vert _{L}=\sup \{\frac{\rho (f(x),f(y))}{d(x,y)}\mid
x,y\in T,x\neq y\}\text{.}
\end{equation*}%
In case $\left\Vert f\right\Vert _{L}<\infty $, we say that $f$ is
lipschitzian. In this case, we have $\rho (f(x),f(y))\leq \left\Vert
f\right\Vert _{L}d(x,y)$ for any $x$ and $y$ in $T$. The set of all
lipschitzian functions $f:T\rightarrow X$ will be denoted by $Lip(T,X)$. In
case $X=T$, we write $Lip(T)$ instead of $Lip(T,X)$. In the particular case
when $X$ is a normed space, it follows that $Lip(T,X)$ is a vector space
seminormed with the seminorm $f\rightarrow \left\Vert f\right\Vert _{L}$. In
the particular case when $(T,d)$ is a compact metric space and $X$ is a
normed space, it follows that $Lip(T,X)\subset C(T,X)\subset B(T,X)$ and $%
Lip(T,X)$ is a normed space with the norm $f\rightarrow \left\Vert
f\right\Vert _{BL}\overset{def}{=}\left\Vert f\right\Vert _{\infty
}+\left\Vert f\right\Vert _{L}$. In the same context, we introduce the sets $%
L_{1}(X)=\{f\in Lip(T,X)\mid \left\Vert f\right\Vert _{L}\leq 1\}$ and $%
BL_{1}(X)=\{f\in Lip(T,X)\mid \left\Vert f\right\Vert _{BL}\leq 1\}$
(clearly $BL_{1}(X)\subset L_{1}(X)$).

A function $f\in Lip(T)$ with $\left\Vert f\right\Vert _{L}<1$ is called a
contraction (with contraction factor $\left\Vert f\right\Vert _{L}$). The
fundamental theorem of the fixed point theory is:

\textit{The Contraction Principle (Banach-Caccioppoli-Picard)}. Assume that $%
(T,d)$ is a complete metric space and $f:T\rightarrow T$ is a contraction.
Then $f$ has an unique fixed point $x^{\ast }\in X$, i.e. $f(x^{\ast
})=x^{\ast }$.

We use standard facts concerning general measure and integral theory. Let us
mention only the fact that, if $\mu $ is an arbitrary positive measure, the
space $L^{2}(\mu )$ with standard norm $\left\Vert .\right\Vert _{2}$ is a
Hilbert space, the scalar product of two elements $\overset{\sim }{f}$ and $%
\overset{\sim }{g}$ in $L^{2}(\mu )$ being $<\overset{\sim }{f},\overset{%
\sim }{g}>=\dint f\overline{g}d\mu $ where $f\in \overset{\sim }{f}$ and $%
g\in \overset{\sim }{g}$ are arbitrary representatives.

Passing to vector measures, we consider an arbitrary non empty set $T$, an
arbitrary $\sigma $-algebra of sets $\mathcal{B\subset P(T)}$ and an
arbitrary Banach space $X$. For an arbitrary $\sigma $-additive measure $\mu
:\mathcal{B}\rightarrow X$, the total variation $\left\vert \mu \right\vert
(T)$ is defined via 
\begin{equation*}
\left\vert \mu \right\vert (T)=\sup \{\underset{i\in I}{\dsum }\left\Vert
\mu (A_{i})\right\Vert \mid (A_{i})_{i\in I}\in Part(T)\}\text{.}
\end{equation*}%
Here $Part(T)$ is the set of all partitions of $T$ (recall that a partition
of $T$ is a finite family $(A_{i})_{i\in I}$ of disjoint sets $A_{i}\in 
\mathcal{B}$ such that $\underset{i\in I}{\cup }A_{i}=T$). Let us introduce%
\begin{equation*}
cabv(\mathcal{B},X)=\{\mu :\mathcal{B}\rightarrow X\mid \mu \text{ is }%
\sigma \text{-additive and }\left\vert \mu \right\vert (T)<\infty \}
\end{equation*}%
which becomes a Banach space, when equipped with the variational norm $\mu
\rightarrow \left\Vert \mu \right\Vert =\left\vert \mu \right\vert (T)$.
Notice that, if $(\mu _{n})_{n}$ is a sequence in $cabv(\mathcal{B},X)$ and $%
\mu \in cabv(\mathcal{B},X)$ is such that $\mu _{n}\underset{n}{\rightarrow }%
\mu $ in $cabv(\mathcal{B},X)$, then $\mu _{n}\underset{n}{\rightarrow }\mu $
pointwise (i.e. $\mu _{n}\underset{n}{(B)\rightarrow }\mu (B)$ for any $B\in 
\mathcal{B}$). For any $0<a<\infty $, write 
\begin{equation*}
B_{a}(X)=\{\mu \in cabv(\mathcal{B},X)\mid \left\Vert \mu \right\Vert \leq
a\}\text{.}
\end{equation*}%
In a similar way, one computes $\left\vert \mu \right\vert (B)$=the
variation of $\mu $ over $B$ for any $B\in \mathcal{B}$.

In the present paper, we work in the particular case when $(T,d)$ is a
compact metric space and we shall write $\mathcal{B}\overset{\text{def}}{=}$%
the Borel sets of $T$. Also, we shall write only $cabv(X)$ instead of $cabv(%
\mathcal{B},X)$.

We continue introducing the basic facts from our previous papers [5] and
[6]. Again $(T,d)$ is a compact metric space and $X$ is a Hilbert space.

A function $f$ of the form $f=\overset{m}{\underset{i=1}{\sum }}\varphi
_{A_{i}}x_{i}$, with $(A_{i})_{1\leq i\leq n}$ forming a partition of $T$
and all $x_{i}\in X$, is called simple. A function $g:T\rightarrow X$ having
the property that there exists a sequence $(f_{n})_{n}$ of simple functions
such that $f_{n}\overset{\text{u}}{\underset{n}{\rightarrow }}f$ (i.e. $%
(f_{n})_{n}$ converges uniformly to $f$) is called totally measurable.

The vector space of totally measurable functions will be denote by $TM(X)$.
We have the inclusion $C(X)\subset TM(X)\subset B(X)$.

In connection with this inclusion, we give some more details. Namely, for a
given $f\in C(X)$, we shall construct the canonical sequence $(f_{m})_{m}$
of simple functions such that $f_{m}\overset{\text{u}}{\underset{m}{%
\rightarrow }}f$. For the compact set $f(T)$, let us fix $m\in \mathbb{N}%
^{\ast }$ and find $y_{1}^{m}=f(t_{1}^{m})$, $y_{2}^{m}=f(t_{2}^{m})$, ...., 
$y_{k(m)}^{m}=f(t_{k(m)}^{m})$ in $f(T)$ such that $f(T)\subset \overset{k(m)%
}{\underset{i=1}{\cup }}B(y_{i}^{m},\frac{1}{m})$ ($B(x,r)$ is the open ball
of centre $x$ and radius $r$). Then $t_{i}^{m}\in
A_{i}^{m}=f^{-1}(B(y_{i}^{m},\frac{1}{m}))\in \mathcal{B}$ with $\overset{%
k(m)}{\underset{i=1}{\cup }}A_{i}=T$. Retaining only the non empty sets, we
obtain the partition $(B_{1}^{m},B_{2}^{m},...,B_{k(m)}^{m})$ of $T$ given
by $B_{1}^{m}=A_{1}^{m}$, ..., $B_{p}^{m}=A_{p}^{m}\smallsetminus \overset{%
p-1}{\underset{i=1}{\cup }}A_{i}^{m}$, abusive notation. Finally define $%
f_{m}=\overset{k(m)}{\underset{i=1}{\dsum }}\varphi _{B_{i}^{m}}z_{i}^{m}$,
where $z_{i}^{m}$ is arbitrarily taken in each $f(B_{i}^{m})$.

For any simple function $f=\overset{m}{\underset{i=1}{\sum }}\varphi
_{A_{i}}x_{i}$ and any $\mu \in cabv(X)$, the integral of $f$\ with respect
to $\mu $ is defined via%
\begin{equation*}
\int fd\mu \overset{def}{=}\overset{m}{\underset{i=1}{\sum }}<x_{i},\mu
(A_{i})>\text{.}
\end{equation*}%
Then, taking an arbitrary $f\in TM(X)$, we extend the previous definition.
Namely, the integral of $f$ with respect to $\mu $ is (coherent definition)%
\begin{equation*}
\int fd\mu =\underset{m}{\lim }\int f_{m}d\mu \text{,}
\end{equation*}%
where $(f_{m})_{m}$ is a sequence of simple functions such that $f_{m}%
\underset{m}{\overset{\text{u}}{\rightarrow }}f$. So our integral is
uniform. It is sesquilinear, because the function $(f,\mu )\mapsto \int
fd\mu $ is linear in $f$ and antilinear in $\mu $ (when work with $K=\mathbb{%
C}$); for $K=\mathbb{R}$ we have bilinearity. Because of the inequality%
\begin{equation*}
\left\vert \int fd\mu \right\vert \leq \left\Vert \mu \right\Vert \left\Vert
f\right\Vert _{\infty }
\end{equation*}%
we see that the aforementioned function of $(f,\mu )$ is continuous for $%
f\in TM(X)$ normed with $\left\Vert .\right\Vert _{\infty }$ and $\mu \in
cabv(X)$ normed with the variational norm.

An important interpretation of the integral just introduced is the fact that
we have an isometric and antilinear isomorphism (bijection) $%
H:cabv(X)\rightarrow C(X)^{^{\prime }}$ which permits the identification $%
cabv(X)\equiv C(X)^{^{\prime }}$. Namely, $H$ acts via $H(\mu )=V_{\mu }$,
where $V_{\mu }(f)=\int fd\mu $ for any $\mu \in cabv\left( X\right) $ and
any $f\in C(X)$. We use the Riesz-Fr\'{e}chet representation theorem
(antilinear identification $X\equiv X^{^{\prime }}$) and the Dinculeanu
theorem (linear identification $C(X)^{^{\prime }}\equiv cabv(X^{^{\prime }})$%
, see [7]).

Using this integral, we introduce on $cabv(X)$ and on some of its subspaces
new norms (weaker than the variational norm).

For any\textit{\ }$\mu \in cabv(X)$, the Monge-Kantorovich norm of $\mu $ is
defined via

\begin{equation*}
\left\Vert \mu \right\Vert _{MK}=\sup \{\left\vert \int fd\mu \right\vert
\mid f\in BL_{1}(X)\}\text{.}
\end{equation*}%
and we get the (generally incomplete) normed space $(cabv(X),\left\Vert
.\right\Vert _{MK})$. For any $\mu \in cabv(X)$ and any $f\in Lip(T,X)$, one
has%
\begin{equation*}
\left\Vert \mu \right\Vert _{MK}\leq \left\Vert \mu \right\Vert \text{ and }%
\left\vert \int fd\mu \right\vert \leq \left\Vert \mu \right\Vert
_{MK}\left\Vert f\right\Vert _{BL}\text{.}
\end{equation*}

For any $v\in X$, let us define 
\begin{equation*}
cabv(X,v)=\{\mu \in cabv(X)\mid \mu (T)=v\}\text{.}
\end{equation*}%
It is clear that $cabv(X,0)$ is a vector subspace of $cabv(X)$ and $\delta
_{t}v\in cabv(X,v)$ for any $t\in T$. It follows that, if $0<a<\infty $ and $%
v\in X$ is such that $\left\Vert v\right\Vert \leq a$, then%
\begin{equation*}
B_{a}(X,v)\overset{def}{=}B_{a}(X)\cap cabv(X,v)
\end{equation*}%
is not empty, because $\delta _{t}v\in B_{a}(X,v)$ for any $t\in T$.

For any $\mu \in cabv(X,0)$, the modified Monge-Kantorovich norm of $\mu $
is defined via 
\begin{equation*}
\left\Vert \mu \right\Vert _{MK}^{\ast }\overset{def}{=}\sup \{\left\vert
\int fd\mu \right\vert \mid f\in L_{1}(X)\}
\end{equation*}%
and we get the (generally incomplete) normed space $(cabv(X,0),\left\Vert
.\right\Vert _{MK}^{\ast })$. For any $\mu \in cabv(X,0)$ and any $f\in
Lip(T,X)$, one has%
\begin{equation*}
\left\vert \int fd\mu \right\vert \leq \left\Vert \mu \right\Vert
_{MK}^{\ast }\left\Vert f\right\Vert _{L}
\end{equation*}

\begin{equation*}
\left\Vert \mu \right\Vert _{MK}\leq \left\Vert \mu \right\Vert _{MK}^{\ast
}\leq \left\Vert \mu \right\Vert diam(T)
\end{equation*}%
\begin{equation*}
\left\Vert \mu \right\Vert _{MK}\leq \left\Vert \mu \right\Vert _{MK}^{\ast
}\leq \left\Vert \mu \right\Vert _{MK}(diam(T)+1)\text{,}
\end{equation*}%
where, as usual, $diam(T)=\sup \{d(x,y)\mid x,y\in T\}$.

Using the aforementioned identification $cabv(X)\equiv C(X)^{^{\prime }}$,
we have the following results, valid for $0<a<\infty $, $n\in \mathbb{N}%
^{\ast }$ and $v\in K^{n}$ with $\left\Vert v\right\Vert \leq a$:

The set $B_{a}(K^{n})$, equipped with the metric $d_{MK}$ given via $%
d_{MK}(\mu ,\nu )=\left\Vert \mu -\nu \right\Vert _{MK}$ and the non empty
set $B_{a}(K^{n},v)$, equipped with the metric $d_{MK}$ or with the
equivalent metric $d_{MK}^{\ast }$ given via $d_{MK}^{\ast }(\mu ,\nu
)=\left\Vert \mu -\nu \right\Vert _{MK}^{\ast }$, are compact metric spaces,
their topology being exactly the topology induced by the weak$^{\ast }$
topology.

In the particular case\textit{\ }$K=\mathbb{R}$\textit{, }$n=1$\ and $a=v=1$%
, the set $B_{1}^{+}(\mathbb{R},1)=B_{1}(\mathbb{R},1)\cap \{\mu :\mathcal{B}%
\rightarrow \mathbb{R})\mid \mu \geq 0\}$ = the probabilities on $\mathcal{B}
$, is weak$^{\ast }$ closed, hence compact for the weak$^{\ast }$ topology
generated by $d_{MK}$ or by $d_{MK}^{\ast }$.

For general topology, see [13]. For functional analysis, see [8]. For
general measure theory, see [11]. For vector measures and integration, see
[7].

\bigskip

\textbf{2. Fractal (Invariant) Vector Measures. The Finite Case}

\bigskip

2.1 \textbf{Framework of the Paragraph}

\bigskip

As previously, we consider a compact metric space $(T,d)$ with Borel sets $%
\mathcal{B}$ and a Hilbert space $X$.

Let $M\in \mathbb{N}$, $M\geq 1$ and $\omega _{i}\in Lip(T)$ with Lipschitz
constants $r_{i}=\left\Vert \omega _{i}\right\Vert _{L}$, $i=1,2,...,M$. In
case all $\omega _{i}$ are contractions, we have $r_{i}<1$, $i=1,2,...,M$.
We say that $(\omega _{1},\omega _{2},...,\omega _{M})$ is an iterated
function system.

Recall that, for any continuous function $h:T\rightarrow T$ and any measure $%
\mu \in cabv(X)$, we can consider the transported measure $h(\mu )\in
cabv(X) $ acting via%
\begin{equation*}
h(\mu )(B)=\mu (h^{-1}(B))
\end{equation*}%
for any $B\in \mathcal{B}$ (we have $\left\Vert h(\mu )\right\Vert \leq
\left\Vert \mu \right\Vert $).

We consider also $R_{i}\in \mathcal{L}(X)$, $i=1,2,...,M$, which together
with $\omega _{i}$, generate the Markov-type operator $H:cabv(X)\rightarrow
cabv(X)$ given via%
\begin{equation*}
H(\mu )=\overset{M}{\underset{i=1}{\sum }}R_{i}\circ \omega _{i}(\mu )
\end{equation*}%
for any $\mu \in cabv(X)$.

Namely $H\in \mathcal{L}(cabv(X))$ with $\left\vert H(\mu )\right\vert
(T)=\left\Vert H(\mu )\right\Vert \leq (\overset{M}{\underset{i=1}{\sum }}%
\left\Vert R_{i}\right\Vert _{o})\left\Vert \mu \right\Vert $, i.e.%
\begin{equation*}
\left\Vert H\right\Vert _{o}\leq \overset{M}{\underset{i=1}{\sum }}%
\left\Vert R_{i}\right\Vert _{o}\text{.}
\end{equation*}

The last assertion follows easily considering a partition $(A_{j})_{1\leq
j\leq n}$ of $T$ and noticing that 
\begin{equation*}
\overset{n}{\underset{j=1}{\sum }}\left\Vert H(\mu )(A_{j})\right\Vert \leq 
\overset{M}{\underset{i=1}{\sum }}\overset{n}{\underset{j=1}{\sum }}%
\left\Vert R_{i}(\mu (\omega _{i}^{-1}(A_{j})))\right\Vert \leq \overset{M}{%
\underset{i=1}{\sum }}\left\Vert R_{i}\right\Vert _{o}\overset{n}{\underset{%
j=1}{\sum }}\left\Vert \mu (\omega _{i}^{-1}(A_{j}))\right\Vert \text{.}
\end{equation*}

\bigskip

\textbf{Lemma 2.1.1. }\textit{Let\ }$f\in L_{1}(X)$\textit{. Define }$g=%
\overset{M}{\underset{i=1}{\sum }}R_{i}^{\ast }\circ f\circ \omega _{i}$%
\textit{. Then }$g\in Lip(T,X)$\textit{\ and we have}%
\begin{equation*}
\left\Vert g\right\Vert _{L}\leq \overset{M}{\underset{i=1}{\sum }}%
\left\Vert R_{i}\right\Vert _{o}r_{i}\text{.}
\end{equation*}

\textit{Proof}. For any $x,y$ in $T$: $\left\Vert g(x)-g(y)\right\Vert \leq 
\overset{M}{\underset{i=1}{\sum }}\left\Vert R_{i}^{\ast }\right\Vert
_{o}\left\Vert f(\omega _{i}(x))-f(\omega _{i}(y))\right\Vert \leq \overset{M%
}{\underset{i=1}{\sum }}\left\Vert R_{i}\right\Vert _{o}\left\Vert \omega
_{i}(x)-\omega _{i}(y)\right\Vert \leq (\overset{M}{\underset{i=1}{\sum }}%
\left\Vert R_{i}\right\Vert _{o}r_{i})d(x,y)$. $\square $

\bigskip

\textbf{Theorem 2.1.2}. \textit{(Change of Variable Formula). For any }$f\in
C(X)$\textit{\ and any }$\mu \in cabv(X)$\textit{, one has}%
\begin{equation*}
\dint fdH(\mu )=\dint gd\mu \text{,}
\end{equation*}%
\textit{where }$g=\overset{M}{\underset{i=1}{\sum }}R_{i}^{\ast }\circ
f\circ \omega _{i}$\textit{.}

\textit{Proof}. Using additivity, it will be sufficient to prove that, for
any $R\in \mathcal{L}(X)$, any continuous $\omega :T\rightarrow T$ and any
continuous $f:T\rightarrow X$, one has%
\begin{equation}
\dint fdH(R)(\mu )=\dint gd\mu \text{,}  \tag{2.1.1}
\end{equation}%
where $H(R)(\mu )\in cabv(X)$ acts via $H(R)(\mu )=R\circ \omega (\mu )$ and 
$g=R^{\ast }\circ f\circ \omega $. Let us construct the canonical sequence $%
(f_{m})_{m}$ for $f$: $f_{m}=\overset{k(m)}{\underset{i=1}{\dsum }}\varphi
_{B_{i}^{m}}f(t_{i}^{m})$, $t_{i}^{m}\in B_{i}^{m}$ and we have $\left\Vert
f_{m}\right\Vert _{\infty }\leq \left\Vert f\right\Vert _{\infty }$ and $%
(C_{i}^{m})_{1\leq i\leq k(m)}$ is a partition of $T$, where $%
C_{i}^{m}=\omega ^{-1}(B_{i}^{m})$.

Take $v_{i}^{m}\in C_{i}^{m}$ with $\omega (v_{i}^{m})=t_{i}^{m}$ and
compute:%
\begin{equation*}
\dint f_{m}dH(R)(\mu )=\overset{k(m)}{\underset{i=1}{\dsum }}%
<f(t_{i}^{m}),H(R)(\mu )(B_{i}^{m})>=
\end{equation*}%
\begin{equation*}
\overset{k(m)}{\underset{i=1}{\dsum }}<f(t_{i}^{m}),R(\mu (\omega
^{-1}(B_{i}^{m})))>=\overset{k(m)}{\underset{i=1}{\dsum }}%
<f(t_{i}^{m}),R(\mu (C_{i}^{m}))>=
\end{equation*}%
\begin{equation*}
=\overset{k(m)}{\underset{i=1}{\dsum }}<(R^{\ast }\circ f)(\omega
(v_{i}^{m})),\mu (C_{i}^{m})>=\overset{k(m)}{\underset{i=1}{\dsum }}%
<(R^{\ast }\circ f\circ \omega )(v_{i}^{m}),\mu (C_{i}^{m})>\text{.}
\end{equation*}

Introducing the simple function $g_{m}=\overset{k(m)}{\underset{i=1}{\dsum }}%
\varphi _{C_{i}^{m}}(R^{\ast }\circ f\circ \omega )(v_{i}^{m})$ we got the
formula%
\begin{equation}
\dint f_{m}dH(R)(\mu )=\dint g_{m}d\mu \text{,}  \tag{2.1.2}
\end{equation}

Now we shall prove that $g_{m}\underset{m}{\overset{\text{u}}{\rightarrow }}%
g $.

Indeed, for any $t\in T$, there exists an unique $i=1,2,...,k(m)$ such that $%
t\in C_{i}^{m}$ and this implies%
\begin{equation*}
\left\Vert g_{m}(t)-g(t)\right\Vert =\left\Vert R^{\ast }(f(\omega
(v_{i}^{m})))-R^{\ast }(f(\omega (t)))\right\Vert =
\end{equation*}%
\begin{equation*}
=\left\Vert R^{\ast }(f(\omega (v_{i}^{m}))-f(\omega (t)))\right\Vert \leq
\left\Vert R\right\Vert _{o}\left\Vert f(\omega (v_{i}^{m}))-f(\omega
(t))\right\Vert \text{.}
\end{equation*}

We have $\omega (v_{i}^{m})=t_{i}^{m}\in B_{i}^{m}$, $\omega (t)\in
B_{i}^{m} $ and $B_{i}^{m}\subset f^{-1}(B(y_{i}^{m},\frac{1}{m}))$ (see the
construction of the canonical sequence $(f_{m})_{m}$), hence $%
f(t_{i}^{m})=f(\omega (v_{i}^{m}))$ and $f(\omega (t))$ are in $B(y_{i}^{m},%
\frac{1}{m})$ and this implies that $\left\Vert f(\omega
(v_{i}^{m}))-f(\omega (t))\right\Vert \leq \frac{2}{m}$, leading to $%
\left\Vert g_{m}(t)-g(t)\right\Vert \leq \left\Vert R\right\Vert _{o}\frac{2%
}{m}$. So $g_{m}\underset{m}{\overset{\text{u}}{\rightarrow }}g$.

Because $\dint fdH(R)(\mu )=\underset{m}{\lim }\dint f_{m}dH(R)(\mu )$, it
follows from $(2.1.2)$ that $(2.1.1)$ is true. $\square $

\bigskip

At the beginning of the chapter, we have seen that, considering on $cabv(X)$
the usual variational norm, the operator $H$ is continuous. In the sequel,
we shall consider on $cabv(X)$ the Monge-Kantorovich norm and we shall see
that $H$ acts as a continuous operator in this context too.

\bigskip

\textbf{Theorem 2.1.3.} \textit{We have} $H\in \mathcal{L}%
(cabv(X),\left\Vert .\right\Vert _{MK})$\textit{. Namely, one has in this
case}%
\begin{equation*}
\left\Vert H\right\Vert _{o}\leq \overset{M}{\underset{i=1}{\dsum }}%
\left\Vert R_{i}\right\Vert _{o}(1+r_{i})\text{.}
\end{equation*}

\textit{Proof}. We take an arbitrary $\mu \in cabv(X)$ and we must prove
that 
\begin{equation}
\left\Vert H(\mu )\right\Vert _{MK}\leq (\overset{M}{\underset{i=1}{\dsum }}%
\left\Vert R_{i}\right\Vert _{o}(1+r_{i}))\left\Vert \mu \right\Vert _{MK}%
\text{.}  \tag{2.1.3}
\end{equation}

Indeed, take an arbitrary $f\in BL_{1}(X)$. Considering the canonical $g=%
\overset{M}{\underset{i=1}{\dsum }}R_{i}^{\ast }\circ f\circ \omega _{i}$,
we see that $\left\Vert g\right\Vert _{\infty }\leq \overset{M}{\underset{i=1%
}{\dsum }}\left\Vert R_{i}\right\Vert _{o}$, because $\left\Vert
g(t)\right\Vert \leq \overset{M}{\underset{i=1}{\dsum }}\left\Vert
R_{i}^{\ast }\right\Vert _{o}\left\Vert f(\omega _{i}(t))\right\Vert $%
\linebreak $\leq \overset{M}{\underset{i=1}{\dsum }}\left\Vert R_{i}^{\ast
}\right\Vert _{o}$, for any $t\in T$. Also, we know (Lemma 2.1.1) that $%
\left\Vert g\right\Vert _{L}\leq \overset{M}{\underset{i=1}{\dsum }}%
\left\Vert R_{i}\right\Vert _{o}r_{i}$. Hence, $\left\Vert g\right\Vert
_{BL}=\left\Vert g\right\Vert _{\infty }+\left\Vert g\right\Vert _{L}\leq 
\overset{M}{\underset{i=1}{\dsum }}\left\Vert R_{i}\right\Vert _{o}(1+r_{i})$
and this implies (use Theorem 2.1.2) 
\begin{equation*}
\left\vert \dint fdH(\mu )\right\vert =\left\vert \dint gd\mu \right\vert
\leq \left\Vert g\right\Vert _{BL}\left\Vert \mu \right\Vert _{MK}\leq (%
\overset{M}{\underset{i=1}{\dsum }}\left\Vert R_{i}\right\Vert
_{o}(1+r_{i}))\left\Vert \mu \right\Vert _{MK}\text{.}
\end{equation*}%
Passing to supremum according to $f\in BL_{1}(X)$, we get $(2.1.3)$. $%
\square $

\bigskip

Working on the subspace $cabv(X,0)$ of $cabv(X)$, equipped with the modified
Monge-Kantorovich norm $\left\Vert .\right\Vert _{MK}^{\ast }$, we can
consider the restriction of $H$ which is again continuous, as the following
result shows.

\bigskip

\textbf{Theorem 2.1.4.} \textit{For any} $\mu \in cabv(X,0)$\textit{, one
has }$H(\mu )\in cabv(X,0)$\textit{. Hence, one can define }$%
H_{0}:cabv(X,0)\rightarrow cabv(X,0)$\textit{, via }$H_{0}(\mu )=H(\mu )$%
\textit{\ and we have }$H_{0}\in \mathcal{L}(cabv(X,0),\left\Vert
.\right\Vert _{MK}^{\ast })$\textit{. Moreover, in this context: }%
\begin{equation*}
\left\Vert H_{0}\right\Vert _{o}\leq \overset{M}{\underset{i=1}{\dsum }}%
\left\Vert R_{i}\right\Vert _{o}r_{i}\text{.}
\end{equation*}

\textit{Proof}. a) If $\mu \in cabv(X,0)$, we have $H(\mu )(T)=\overset{M}{%
\underset{i=1}{\sum }}R_{i}(\mu (\omega _{i}^{-1}(T)))=\overset{M}{\underset{%
i=1}{\sum }}R_{i}(\mu (T))=0$.

b) Take arbitrarily $\mu \in cabv(X,0)$. We must show that%
\begin{equation}
\left\Vert H(\mu )\right\Vert _{MK}^{\ast }\leq \overset{M}{\underset{i=1}{%
\dsum }}\left\Vert R_{i}\right\Vert _{o}r_{i})\left\Vert \mu \right\Vert
_{MK}^{\ast }\text{.}  \tag{2.1.4}
\end{equation}

Indeed, take arbitrary $f\in L_{1}(X)$ and construct the canonical\linebreak 
$g=\overset{M}{\underset{i=1}{\dsum }}R_{i}^{\ast }\circ f\circ \omega _{i}$
with $\left\Vert g\right\Vert _{L}\leq \overset{M}{\underset{i=1}{\dsum }}%
\left\Vert R_{i}\right\Vert _{o}r_{i}$ (according to Lemma 2.1.1). Then,
using Theorem 2.1.2:%
\begin{equation*}
\left\vert \dint fdH(\mu )\right\vert =\left\vert \dint fgd\mu \right\vert
\leq \left\Vert g\right\Vert _{L}\left\Vert \mu \right\Vert _{MK}^{\ast
}\leq (\overset{M}{\underset{i=1}{\dsum }}\left\Vert R_{i}\right\Vert
_{o}r_{i}))\left\Vert \mu \right\Vert _{MK}^{\ast }\text{.}
\end{equation*}%
Passing to supremum according to $f\in L_{1}(X)$, we get $(2.1.4)$. $\square 
$

\bigskip

In view of the preceding facts, we shall use the operator $H$ and we shall
introduce two constructions, producing two models, which will be illustrated
further.

\newpage

\textbf{Construction 1} ($H_{1}$ - model)

Let $\emptyset \neq A\subset cabv(X)$ and assume that $H(A)\subseteq A$.
Define $H_{1}:A\rightarrow A$ via $H_{1}(\mu )=H(\mu )$ for each $\mu \in A$.

\bigskip

\textbf{Construction 2} ($H_{2}$ - model)

Let $\emptyset \neq A\subset cabv(X)$ and $\mu ^{0}\in cabv(X)$. Assume that 
$H(A)+\mu ^{0}\overset{def}{=}\{H(\mu )+\mu ^{0}\mid \mu \in A\}\subseteq A$%
. Define $H_{2}:A\rightarrow A$ via $H_{2}(\mu )=H(\mu )+\mu _{0}$ for each $%
\mu \in A$.

\bigskip

Any fixed point $\mu ^{\ast }\in A$ of $H_{i}$, i.e. $H_{i}(\mu ^{\ast
})=\mu ^{\ast }$, will be called a fractal (invariant) measure of $H_{i}$,
or a Hutchinson (self-similar) measure of $H_{i}$, $i=1,2$.

\bigskip

2.2. \textbf{Illustrations of the }$H_{1}$\textbf{\ and }$H_{1}$\textbf{\
Models}

\bigskip

All the concrete illustrations in this paragraph will be done within the
following particular framework: $T=[0,1]$, $M=2$ and $\omega _{1},\omega
_{2} $ are the Cantor contractions: 
\begin{equation*}
\omega _{1}:[0,1]\rightarrow \lbrack 0,1],\omega _{1}(t)=\frac{t}{3}\text{,
with }r_{1}=\frac{1}{3}
\end{equation*}%
\begin{equation*}
\omega _{2}:[0,1]\rightarrow \lbrack 0,1],\omega _{1}(t)=\frac{2}{3}+\frac{t%
}{3}\text{, with }r_{2}=\frac{1}{3}\text{.}
\end{equation*}

Hence for any $\emptyset \neq B\in \mathcal{B}$, one has:%
\begin{equation*}
\omega _{1}^{-1}(B)=3B\cap \lbrack 0,1]=\{3t\mid t\in B\}\cap \lbrack 0,1]
\end{equation*}%
\begin{equation*}
\omega _{2}^{-1}(B)=(3B-2)\cap \lbrack 0,1]=\{3t-2\mid t\in B\}\cap \lbrack
0,1]\text{.}
\end{equation*}

\bigskip

In this subparagraph $\lambda :\mathcal{B\rightarrow }\mathbb{R}_{+}$ is the
Lebesgue measure on $[0,1]$.

\bigskip

A. \textbf{Illustration of the }$H_{1}$\textbf{\ - model}

\bigskip

\textbf{Theorem 2.2.1.} \textit{Consider }$X=K^{n}$\textit{, }$n\in \mathbb{N%
}^{\ast }$\textit{. The hypotheses are:}

\textit{a) }$\overset{M}{\underset{i=1}{\sum }}R_{i}=Id_{K^{n}};$

\textit{b) }$c=\overset{M}{\underset{i=1}{\sum }}\left\Vert R_{i}\right\Vert
_{o}r_{i}<1$\textit{\ (clearly this true if all }$\omega _{i}$\textit{\ are
contractions and }$\overset{M}{\underset{i=1}{\sum }}\left\Vert
R_{i}\right\Vert _{o}=1)$\textit{;}

\textit{c) }$0<a<\infty $\textit{\ and }$\upsilon \in K^{n}$\textit{\ are
such that }$\left\Vert \upsilon \right\Vert \leq a$\textit{;}

\textit{d) }$\emptyset \neq A\subseteq B_{a}(K^{n},v)$\textit{\ is such that 
}$H(A)\subseteq A$\textit{\ and }$A$\textit{\ is weak}$^{\ast }$\textit{\
closed (In the particular case when }$\left\Vert H(\mu )\right\Vert \leq
\left\Vert \mu \right\Vert $\textit{\ for any }$\mu \in cabv(K^{n})$\textit{%
, one can take }$A=B_{a}(K^{n},v)$\textit{. More particular, if }$\overset{M}%
{\underset{i=1}{\sum }}\left\Vert R_{i}\right\Vert _{o}=1$\textit{, it
follows that }$\left\Vert H(\mu )\right\Vert \leq \left\Vert \mu \right\Vert 
$\textit{\ for any }$\mu \in cabv(K^{n})$\textit{).}

\textit{Under these hypotheses, we define }$H_{1}:A\rightarrow A$\textit{\
via }$H_{1}(\mu )=H(\mu )$\textit{\ for any }$\mu \in A$\textit{. It follows
that }$H_{1}$\textit{\ is a contraction with contraction factor }$\leq c$%
\textit{, if }$A$\textit{\ is equipped with the metric }$d_{MK}^{\ast }$%
\textit{\ given via }$d_{MK}^{\ast }(\mu ,\nu )=\left\Vert \mu -\nu
\right\Vert _{MK}^{\ast }$\textit{.}

\textit{There exists an unique fractal (invariant) measure }$\mu ^{\ast }\in
A$\textit{\ of }$H_{1}$\textit{, i.e. }$H_{1}(\mu ^{\ast })=\mu ^{\ast }$%
\textit{.}

\textit{Proof}. We shall prove the general assertion, the particular cases
being discussed at the end.

According to the Preliminary Part, $B_{a}(K^{n},v)$ is a non empty compact
space for the metric $d_{MK}^{\ast }$. Consequently, $A$ is also compact for
this metric, being weak$^{\ast }$ closed, hence closed in $B_{a}(K^{n},v)$
(coincidence of the weak$^{\ast }$ topology with the topology given by $%
d_{MK}^{\ast }$).

Condition a) guarantees that $H(cabv(K^{n},v))\subset cabv(K^{n},v)$: if $%
\mu \in cabv(K^{n},v)$, then $H(\mu )(T)=\overset{M}{\underset{i=1}{\sum }}%
R_{i}(\mu (\omega _{i}^{-1}(T)))=\overset{M}{\underset{i=1}{\sum }}R_{i}(\mu
(T))=\overset{M}{\underset{i=1}{\sum }}R_{i}(v)=v$.

Now we prove that $H_{1}$ is a contraction with contraction factor $\leq C$.
To this end, take $\mu $ and $\nu $ in $A$. Then $\mu -\nu \in cabv(K^{n},0)$%
, hence $H(\mu -\nu )\in cabv(K^{n},0)$ and%
\begin{equation*}
d_{MK}^{\ast }(H_{1}(\mu ),H_{1}(\nu ))=\left\Vert H_{1}(\mu )-H_{1}(\nu
)\right\Vert _{MK}^{\ast }=\left\Vert H(\mu -\nu )\right\Vert _{MK}^{\ast
}\leq
\end{equation*}%
\begin{equation*}
\leq \overset{M}{\underset{i=1}{(\sum }}\left\Vert R_{i}\right\Vert
_{o}r_{i})\left\Vert \mu -\nu \right\Vert _{MK}^{\ast }=cd_{MK}^{\ast }(\mu
,\nu )
\end{equation*}%
according to Theorem 2.1.4.

The existence and uniqueness of $\mu ^{\ast }$ follow from the contraction
principle.

Concerning the particular cases, we see first that, in case $\left\Vert
H(\mu )\right\Vert \leq \left\Vert \mu \right\Vert $ for any $\mu \in
cabv(K^{n})$, we have $H_{1}(B_{a}(K^{n},v))\subset
B_{a}(K^{n},v)=B_{a}(K^{n})\cap cabv(X,v)$.

Finally, if $\overset{M}{\underset{i=1}{\sum }}\left\Vert R_{i}\right\Vert
_{o}=1$, take $\mu \in cabv(K^{n})$ and use the evaluation at the beginning
of the paragraph%
\begin{equation*}
\left\Vert H(\mu )\right\Vert \leq \overset{M}{\underset{i=1}{\dsum }}%
\left\Vert R_{i}\right\Vert _{o})\left\Vert \mu \right\Vert =\left\Vert \mu
\right\Vert \text{. }\square
\end{equation*}

\bigskip

\textbf{Remarks}

\bigskip

1. Condition a) implies that $1=\left\Vert Id_{K^{n}}\right\Vert _{o}\leq 
\overset{M}{\underset{i=1}{\sum }}\left\Vert R_{i}\right\Vert _{o}$, hence
the condition $\overset{M}{\underset{i=1}{\sum }}\left\Vert R_{i}\right\Vert
_{o}=1$\textit{\ }is extremal.

There exist situations when all the particular conditions are fulfilled, as
we can see in the following remark.

\bigskip

2. The classical model, producing the fractal (invariant) probability is a
particular case of Theorem 2.2.1 where all the particular conditions are
fulfilled.

Namely one takes $n=1$ (hence $X=K$), $R_{i}\in \mathcal{L}(K)$ given via $%
R_{i}(t)=p_{i}t$, where all $p_{i}>0$, $i=1,2,...,M$ and $\overset{M}{%
\underset{i=1}{\sum }}p_{i}=1$, hence $H(\mu )=\overset{M}{\underset{i=1}{%
\sum }}p_{i}\omega _{i}(\mu )$. Also, one takes $a=1$, $v=1$ and $A=\{\mu
\in B_{1}(K,1)\mid \mu \geq 0\}=$the probabilities $\mu :\mathcal{%
B\rightarrow }[0,1]$. Then $A$ is weak$^{\ast }$ closed (see the Preliminary
Part) and, for any contraction $\omega _{i}:T\rightarrow T$, $i=1,2,...,M$
one has: $\overset{M}{\underset{i=1}{\sum }}R_{i}=Id_{K}$, $\overset{M}{%
\underset{i=1}{\sum }}\left\Vert R_{i}\right\Vert _{o}=\overset{M}{\underset{%
i=1}{\sum }}p_{i}=1$ (hence $C<1$). We find an unique probability $\mu
^{\ast }:\mathcal{B\rightarrow }[0,1]$ (the fractal invariant measure)
having the property $\mu ^{\ast }=\overset{M}{\underset{i=1}{\sum }}%
p_{i}\omega _{i}(\mu ^{\ast })$.

\bigskip

\textbf{Concrete Illustrations}

\bigskip

Take $n=2$ (hence $X=K^{2}$) and $R_{1}$, $R_{2}$ in $\mathcal{L}(K^{2})$
such that 
\begin{equation*}
R_{1}\equiv (%
\begin{array}{cc}
\alpha & 0 \\ 
0 & \alpha%
\end{array}%
)\text{, }R_{2}\equiv (%
\begin{array}{cc}
1-\alpha & 0 \\ 
0 & 1-\alpha%
\end{array}%
)\text{,}
\end{equation*}
where $0<\alpha <1$. It follows that $R_{1}+R_{2}=Id_{K^{2}}$, $\left\Vert
R_{1}\right\Vert _{o}=\alpha $, $\left\Vert R_{2}\right\Vert _{o}=1-\alpha $%
, hence $\left\Vert R_{1}\right\Vert _{o}+\left\Vert R_{2}\right\Vert _{o}=1$%
. Also take $a=\sqrt{2}$ and $v=(1,1)$, hence $\left\Vert v\right\Vert =a$.

We get the fractal (invariant) measure $\mu ^{\ast }=(\mu _{1}^{\ast },\mu
_{2}^{\ast })$. Namely, the invariance equation $H_{1}(\mu ^{\ast })=\mu
^{\ast }$, i.e. $R_{1}\circ \omega _{1}(\mu ^{\ast })+R_{2}\circ \omega
_{2}(\mu ^{\ast })=\mu ^{\ast }$ is (for any $B\in \mathcal{B}$):%
\begin{equation*}
R_{1}(\mu ^{\ast }((3B)\cap \lbrack 0,1]))+R_{2}(\mu ^{\ast }(3B-2)\cap
\lbrack 0,1]))=\mu ^{\ast }(B)\text{.}
\end{equation*}

In matricial form%
\begin{equation*}
(%
\begin{array}{cc}
\alpha & 0 \\ 
0 & \alpha%
\end{array}%
)(%
\begin{array}{c}
\mu _{1}^{\ast }((3B)\cap \lbrack 0,1]) \\ 
\mu _{2}^{\ast }((3B)\cap \lbrack 0,1])%
\end{array}%
)+(%
\begin{array}{cc}
1-\alpha & 0 \\ 
0 & 1-\alpha%
\end{array}%
)(%
\begin{array}{c}
\mu _{1}^{\ast }((3B-2)\cap \lbrack 0,1]) \\ 
\mu _{2}^{\ast }((3B-2)\cap \lbrack 0,1])%
\end{array}%
)=
\end{equation*}%
\begin{equation*}
=(%
\begin{array}{c}
\mu _{1}^{\ast }(B) \\ 
\mu _{2}^{\ast }(B)%
\end{array}%
)\text{,}
\end{equation*}%
giving 
\begin{equation*}
\alpha \mu _{i}^{\ast }((3B)\cap \lbrack 0,1])+(1-\alpha )\mu _{i}^{\ast
}((3B-2)\cap \lbrack 0,1])=\mu _{i}^{\ast }(B)\text{,}
\end{equation*}%
$i=1,2$.

Hence $\mu _{1}^{\ast }=\mu _{2}^{\ast }=\mu $, where $\mu :\mathcal{B}%
\rightarrow \lbrack 0,1]$ is the unique fractal (invariant) probability
obtained in the classic model for $p_{1}=\alpha $ and $p_{2}=1-\alpha $.

\bigskip

B. \textbf{First} \textbf{Illustration of the }$H_{2}$\textbf{\ - model}

\bigskip

\textbf{Theorem 2.2.2.} \textit{Consider }$X=K^{n}$\textit{, }$n\in \mathbb{N%
}^{\ast }$\textit{. The hypotheses are:}

\textit{a) }$d=\overset{M}{\underset{i=1}{\sum }}\left\Vert R_{i}\right\Vert
_{o}(1+r_{i})<1$\textit{\ (This true, in particular, if all }$\omega _{i}$%
\textit{\ are contractions and }$\overset{M}{\underset{i=1}{\sum }}%
\left\Vert R_{i}\right\Vert _{o}\leq \frac{1}{2})$\textit{;}

\textit{b) }$0<a<\infty $, $\mu ^{0}\in cabv(K^{n})$, $\emptyset \neq
A\subseteq B_{a}(K^{n})$\textit{\ is weak}$^{\ast }$\textit{\ closed and one
has }$H(\mu )+\mu ^{0}\in A$ \textit{for any} $\mu \in A$ \textit{(In
particular, if }$\left\Vert \mu ^{0}\right\Vert +a(\overset{M}{\underset{i=1}%
{\sum }}\left\Vert R_{i}\right\Vert _{o})\leq a$,\textit{\ then one can take 
}$A=B_{a}(K^{n})$\textit{).}

\textit{Under these hypotheses, we define }$H_{2}:A\rightarrow A$\textit{\
via }$H_{2}(\mu )=H(\mu )+\mu ^{0}$\textit{\ for any }$\mu \in A$\textit{.
It follows that }$H_{2}$\textit{\ is a contraction with contraction factor }$%
\leq d$\textit{, if }$A$\textit{\ is equipped with the metric }$d_{MK}$%
\textit{, given via }$d_{MK}(\mu ,\nu )=\left\Vert \mu -\nu \right\Vert
_{MK} $\textit{.}

\textit{Then:}

\textit{i) If }$\mu ^{0}=0$\textit{, it follows that }$0\in A$\textit{.}

\textit{ii) There exists an unique fractal (invariant) measure }$\mu ^{\ast
}\in A$\textit{\ of }$H_{2}$\textit{, i.e. }$H_{2}(\mu ^{\ast })=\mu ^{\ast
} $\textit{. In case }$\mu ^{0}=0$\textit{, we have }$\mu ^{\ast }=0$.

\textit{Proof}. Again we use the fact that $B_{a}((K^{n})$ is a compact
metric space when being equipped with $d_{MK}$, hence $A$ is in the same
situation, being closed in $B_{a}((K^{n})$.

Let us prove that $H_{2}$ is a contraction with contraction factor $\leq d$.
To this end, take $\mu $ and $\nu $ in $A$. We have%
\begin{equation*}
d_{MK}(H_{2}(\mu ),H_{2}(\nu ))=\left\Vert H_{2}(\mu )-H_{2}(\nu
)\right\Vert _{MK}=\left\Vert H(\mu -\nu )\right\Vert _{MK}\leq
\end{equation*}%
\begin{equation*}
\leq (\overset{M}{\underset{i=1}{\sum }}\left\Vert R_{i}\right\Vert
_{o}(1+r_{i}))\left\Vert \mu -\nu \right\Vert _{MK}=d\left\Vert \mu -\nu
\right\Vert _{MK}=dd_{MK}(\mu ,\nu )
\end{equation*}%
according to Theorem 2.1.3.

The existence and uniqueness of $\mu ^{\ast }$ follow from the contraction
principle.

The particular case concerning point b) is treated as follows. Take $\mu \in
B_{a}(K^{n})$. Then, using the evaluation from the beginning of the
paragraph: $\left\Vert H(\mu )+\mu ^{0}\right\Vert \leq \left\Vert H(\mu
)\right\Vert +\left\Vert \mu ^{0}\right\Vert \leq (\overset{M}{\underset{i=1}%
{\sum }}\left\Vert R_{i}\right\Vert _{o})\left\Vert \mu \right\Vert
+\left\Vert \mu ^{0}\right\Vert \leq a(\overset{M}{\underset{i=1}{\sum }}%
\left\Vert R_{i}\right\Vert _{o})+\left\Vert \mu ^{0}\right\Vert \leq a$,
hence $H(\mu )+\mu ^{0}\in B_{a}(K^{n})$.

The proof finishes with the study of the case when $\mu ^{0}=0$.

First we show that in this case one must have $0\in A$.

Indeed, we have successively: $H(\mu )\in A$, $H(H(\mu ))=H^{2}(\mu )\in A$,
..., $H^{n}(\mu )\overset{def}{=}H(H^{n-1}(\mu ))\in A$ for any $n\in 
\mathbb{N}^{\ast }$. But, according to Theorem 2.1.3: $\left\Vert
H\right\Vert _{o}\leq d$, $\left\Vert H^{2}\right\Vert _{o}\leq d^{2}$, ..., 
$\left\Vert H^{n}\right\Vert _{o}\leq d^{n}\underset{n}{\rightarrow }0$,
hence $\left\Vert H^{n}(\mu )\right\Vert _{MK}\leq \left\Vert
H^{n}\right\Vert _{o}\left\Vert \mu \right\Vert _{MK}\underset{n}{%
\rightarrow }0$. Because $A$ is weak$^{\ast }$ closed, it is closed in the
topology generated by $\left\Vert .\right\Vert _{MK}$ too, consequently $%
\underset{n}{\lim }H^{n}(\mu )=0\in A$.

Because $0\in A$, we can write $H_{2}(0)=0$, so $0$ is a fixed point for $%
H_{2}$ and, due to uniqueness, we have $\mu ^{\ast }=0$. $\square $

\bigskip

\textbf{Concrete Illustrations}

\bigskip

Take $n=2$ (hence $X=K^{2}$) and $\mu ^{0}:\mathcal{B}\rightarrow K^{2}$
acting via $\mu ^{0}(B)=(\frac{1}{4}\lambda (B),\frac{1}{4}\delta _{0}(B))$
for any $B\in \mathcal{B}$. Here $\delta _{0}:\mathcal{B}\rightarrow \mathbb{%
R}_{+}$ is the Dirac measure on $[0,1]$ concentrated at $0$. Take also $%
R_{1} $, $R_{2}$ in $\mathcal{L}(K^{2})$ as follows: $R_{1}=\frac{1}{10}P_{1}
$, $R_{2}=\frac{1}{10}P_{2}$, where $P_{1}$, $P_{2}$ in $\mathcal{L}(K^{2})$
are such that 
\begin{equation*}
P_{1}\equiv (%
\begin{array}{cc}
1 & 0 \\ 
2 & 1%
\end{array}%
)\text{, }P_{2}\equiv (%
\begin{array}{cc}
1 & 0 \\ 
2 & -1%
\end{array}%
)\text{,}
\end{equation*}%
consequently $\left\Vert P_{1}\right\Vert _{o}=\left\Vert P_{2}\right\Vert
_{o}=1+\sqrt{2}$, giving $\left\Vert R_{1}\right\Vert _{o}+\left\Vert
R_{2}\right\Vert _{o}=\frac{1+\sqrt{2}}{5}<\frac{1}{10}$, so $d<1$.

Take $a=1$, hence $\left\Vert \mu ^{0}\right\Vert +a(\overset{2}{\underset{%
i=1}{\sum }}\left\Vert R_{i}\right\Vert _{o})=\frac{1}{2}+\frac{1+\sqrt{2}}{5%
}<1=a$ (We accept that $\left\Vert \mu ^{0}\right\Vert =\frac{1}{2}$ and
this is proved as follows: write $\nu ^{0}=(\lambda ,\delta _{0})$, i.e. $%
\mu ^{0}=\frac{1}{4}\nu ^{0}$. We have: $\left\Vert \nu ^{0}\right\Vert
=\left\vert \nu ^{0}\right\vert ([0,1])=\left\vert \nu ^{0}\right\vert
(\{0\})+\left\vert \nu ^{0}\right\vert ((0,1])$ and $\left\vert \nu
^{0}\right\vert (\{0\})=\left\Vert (\lambda (\{0\}),\delta
_{0}(\{0\}))\right\Vert =\left\Vert (0,1)\right\Vert =1$. For any $%
(0,1]\supset B\in \mathcal{B}$ one has $\nu ^{0}(B)=(\lambda (B),\delta
_{0}(B))=(\lambda (B),0)$, hence $\left\Vert \nu ^{0}(B)\right\Vert =\lambda
(B)$ and this leads to $\left\vert \nu ^{0}\right\vert (B)=\lambda (B)$. We
decide that $\left\Vert \nu ^{0}\right\Vert =1+1=2$, so $\left\Vert \mu
^{0}\right\Vert =\frac{2}{4}=\frac{1}{2}$.).

All the conditions in Theorem 2.1.2 are fulfilled and we obtain the fractal
(invariant) measure $\mu ^{\ast }=(\mu _{1}^{\ast },\mu _{2}^{\ast })\in
cabv(K^{2})$.

The invariance equation $H_{2}(\mu ^{\ast })=\mu ^{\ast }$ is, for any $B\in 
\mathcal{B}$: 
\begin{equation*}
R_{1}(\mu ^{\ast }((3B)\cap \lbrack 0,1]))+R_{2}(\mu ^{\ast }(3B-2)\cap
\lbrack 0,1]))+\mu ^{0}(B)=\mu ^{\ast }(B)\text{.}
\end{equation*}

In matricial form:

\begin{equation*}
(%
\begin{array}{cc}
\frac{1}{10} & 0 \\ 
\frac{2}{10} & \frac{1}{10}%
\end{array}%
)(%
\begin{array}{c}
\mu _{1}^{\ast }((3B)\cap \lbrack 0,1]) \\ 
\mu _{2}^{\ast }((3B)\cap \lbrack 0,1])%
\end{array}%
)+(%
\begin{array}{cc}
\frac{1}{10} & 0 \\ 
\frac{2}{10} & -\frac{1}{10}%
\end{array}%
)(%
\begin{array}{c}
\mu _{1}^{\ast }((3B-2)\cap \lbrack 0,1]) \\ 
\mu _{2}^{\ast }((3B-2)\cap \lbrack 0,1])%
\end{array}%
)+
\end{equation*}%
\begin{equation*}
+(%
\begin{array}{c}
\frac{1}{4}\lambda (B) \\ 
\frac{1}{4}\delta _{0}(B)%
\end{array}%
)=(%
\begin{array}{c}
\mu _{1}^{\ast }(B) \\ 
\mu _{2}^{\ast }(B)%
\end{array}%
)\text{,}
\end{equation*}%
giving

\begin{equation*}
\frac{1}{10}\mu _{1}^{\ast }((3B)\cap \lbrack 0,1])+\frac{1}{10}\mu
_{1}^{\ast }((3B-2)\cap \lbrack 0,1])+\frac{1}{4}\lambda (B)=\mu _{1}^{\ast
}(B)\text{,}
\end{equation*}%
\begin{equation*}
\frac{2}{10}\mu _{1}^{\ast }((3B)\cap \lbrack 0,1])+\frac{1}{10}\mu
_{2}^{\ast }((3B)\cap \lbrack 0,1])+\frac{2}{10}\mu _{1}^{\ast }((3B-2)\cap
\lbrack 0,1])-
\end{equation*}%
\begin{equation*}
-\frac{1}{10}\mu _{2}^{\ast }((3B-2)\cap \lbrack 0,1])+\frac{1}{4}\delta
_{0}(B)=\mu _{2}^{\ast }(B)\text{.}
\end{equation*}

\bigskip

Examples of computation

\bigskip

1. Take $B=[0,1]$, hence $(3B)\cap \lbrack 0,1]=(3B-2)\cap \lbrack
0,1])=[0,1]$. Write $\mu _{1}^{\ast }([0,1])=x$ and $\mu _{2}^{\ast
}([0,1])=y$. Then 
\begin{equation*}
\frac{1}{10}x+\frac{1}{10}x+\frac{1}{4}=x
\end{equation*}%
\begin{equation*}
\frac{2}{10}x+\frac{1}{10}y+\frac{2}{10}x-\frac{1}{10}y+\frac{1}{4}=y\text{,}
\end{equation*}%
giving $x=\mu _{1}^{\ast }([0,1])=\frac{5}{16}$ and $y=\mu _{2}^{\ast
}([0,1])=\frac{3}{8}$.

\bigskip

2. Write $\mu _{i}^{\ast }(\{t\})=\mu _{i}^{\ast }(t)$ and let us compute $%
\mu _{i}^{\ast }(t)$ for some $t\in \lbrack 0,1]$, $i=1,2$.

a) For $B=\{0\}$, we get $\frac{1}{10}\mu _{1}^{\ast }(0)=\mu _{1}^{\ast
}(0) $, hence $\mu _{1}^{\ast }(0)=0$ and $\frac{1}{5}\mu _{1}^{\ast }(0)+%
\frac{1}{10}\mu _{2}^{\ast }(0)+\frac{1}{4}=\mu _{1}^{\ast }(0)$, hence $\mu
_{2}^{\ast }(0)=\frac{5}{18}$.

b) For $B=\{1\}$, we get $\frac{1}{10}\mu _{1}^{\ast }(1)=\mu _{1}^{\ast
}(1) $, hence $\mu _{1}^{\ast }(1)=0$ and $\frac{1}{5}\mu _{1}^{\ast }(1)-%
\frac{1}{10}\mu _{2}^{\ast }(1)=\mu _{2}^{\ast }(1)$, hence $\mu _{2}^{\ast
}(1)=0$.

c) For $B=\{\frac{1}{3}\}$, we get $\frac{1}{10}\mu _{1}^{\ast }(1)=\mu
_{1}^{\ast }(\frac{1}{3})$, hence $\mu _{1}^{\ast }(\frac{1}{3})=0$ and $%
\frac{1}{5}\mu _{1}^{\ast }(1)+\frac{1}{10}\mu _{2}^{\ast }(1)=\mu
_{2}^{\ast }(\frac{1}{3})$, hence $\mu _{2}^{\ast }(\frac{1}{3})=0$.

d) For $B=\{\frac{2}{3}\}$, we get $\frac{1}{10}\mu _{1}^{\ast }(0)=\mu
_{1}^{\ast }(\frac{1}{3})$, hence $\mu _{1}^{\ast }(\frac{1}{3})=0$ and $%
\frac{1}{5}\mu _{1}^{\ast }(0)-\frac{1}{10}\mu _{2}^{\ast }(0)=\mu
_{2}^{\ast }(\frac{2}{3})$, hence $\mu _{2}^{\ast }(\frac{2}{3})=-\frac{1}{36%
}$.

\bigskip

C. \textbf{Second} \textbf{Illustration of the }$H_{2}$\textbf{\ - model}

\bigskip

\textbf{Theorem 2.2.3.} \textit{We work in an arbitrary Hilbert space }$X$ 
\textit{(as a matter of fact, this theorem is valid for any Banach space }$X$%
\textit{) and consider the usual Banach space }$cabv(X)$ \textit{with the
variational norm. Take }$\mu ^{0}\in cabv(X)$.

\textit{The hypotheses are:}

\textit{a) }$e=\overset{M}{\underset{i=1}{\sum }}\left\Vert R_{i}\right\Vert
_{o}<1$\textit{;}

\textit{b)} $\emptyset \neq A\subseteq cabv(X)$\textit{\ is a\ closed set
such that }$H(\mu )+\mu ^{0}\in A$ \textit{for any} $\mu \in A$ \textit{(In
particular, one can take }$A=cabv(X)$\textit{\ or one can take }$A=B_{a}(X)$%
, \textit{if }$0<a<\infty $ and $\left\Vert \mu ^{0}\right\Vert +a(\overset{M%
}{\underset{i=1}{\sum }}\left\Vert R_{i}\right\Vert _{o})\leq a$\textit{).}

\textit{Under these hypotheses, define }$H_{2}:A\rightarrow A$\textit{\ via }%
$H_{2}(\mu )=H(\mu )+\mu ^{0}$\textit{\ for any }$\mu \in A$\textit{. It
follows that }$H_{2}$\textit{\ is a contraction with contraction factor }$%
\leq e$\textit{.}

\textit{Then:}

\textit{i) If }$\mu ^{0}=0$\textit{, then }$0\in A$\textit{.}

\textit{ii) There exists an unique fractal (invariant) measure }$\mu ^{\ast
}\in A$\textit{\ of }$H_{2}$\textit{, i.e. }$H_{2}(\mu ^{\ast })=\mu ^{\ast
} $\textit{. In case }$\mu ^{0}=0$\textit{, we have }$\mu ^{\ast }=0$.

\textit{Proof}. Using the fact that $\left\Vert H\right\Vert _{o}\leq 
\overset{M}{\underset{i=1}{\sum }}\left\Vert R_{i}\right\Vert _{o}<1$ (see
the beginning of the paragraph), it is easily seen that $H_{2}$ is a
contraction with contraction factor $\leq e$. The existence and uniqueness
of $\mu ^{\ast }$ follow from the contraction principle.

In case $\left\Vert \mu ^{0}\right\Vert +a(\overset{M}{\underset{i=1}{\sum }}%
\left\Vert R_{i}\right\Vert _{o})\leq a$, we take an arbitrary $\mu \in
B_{a}(X)$ and obtain $\left\Vert H(\mu )+\mu ^{0}\right\Vert \leq \left\Vert
H(\mu )\right\Vert +\left\Vert \mu ^{0}\right\Vert \leq (\overset{M}{%
\underset{i=1}{\sum }}\left\Vert R_{i}\right\Vert _{o})\left\Vert \mu
\right\Vert +\left\Vert \mu ^{0}\right\Vert \leq a$, showing that $H(\mu
)+\mu ^{0}\in B_{a}(X)$.

The study of the situation when $\mu ^{0}=0$ is similar to the study for the
case of Theorem 2.2.2. $\square $

\bigskip

\textbf{Concrete Illustrations}

\bigskip

We begin with some initial facts.

Any continuous function $F:[0,1]^{2}\rightarrow K$ with $Q\overset{def}{=}%
\sup \{\left\vert F(x,y)\right\vert \mid (x,y)\in \lbrack 0,1]^{2}\}$
generates $R:L^{2}(\lambda )\rightarrow L^{2}(\lambda )$ given via $R(%
\overset{\sim }{f})=\overset{\sim }{g}$, where $g:[0,1]\rightarrow K$ acts
as follows (we work with a representative $f\in \overset{\sim }{f}$):%
\begin{equation*}
g(x)=\overset{1}{\underset{0}{\int }}F(x,y)f(y)d\lambda (y)
\end{equation*}%
and $g$ is continuous (because $F$ is uniformly continuous and, for any $%
x,x_{0}$ in $[0,1]$ we have $\left\vert g(x)-g(x_{0})\right\vert \leq \int
\left\vert F(x,y)-F(x_{0},y)\right\vert \left\vert f(y)\right\vert d\lambda
(y)$).

We also have $\left\Vert g\right\Vert _{2}\leq Q\left\Vert f\right\Vert _{2}$
because, if $x\in \lbrack 0,1]$, one has $\left\vert g(x)\right\vert
^{2}\leq (\int \left\vert f(y)\right\vert d\lambda (y))^{2}Q^{2}\leq
Q^{2}(\left\Vert f\right\Vert _{2}\left\Vert 1\right\Vert
_{2})^{2}=Q^{2}\left\Vert f\right\Vert _{2}^{2}$.

Hence $\left\Vert R(\overset{\sim }{f})\right\Vert _{2}\leq Q\left\Vert 
\overset{\sim }{f}\right\Vert _{2}$ and $R$ is continuous with $\left\Vert
R\right\Vert _{o}\leq Q$.

Now we shall introduce our example.

We shall take a number $a>0$, $X=L^{2}(\lambda )$ and $F_{i}:[0,1]^{2}%
\rightarrow K$, continuous, with $Q_{i}\overset{def}{=}\sup \{\left\vert
F_{i}(x,y)\right\vert \mid (x,y)\in \lbrack 0,1]^{2}\}$ and we shall assume
that $Q_{i}\leq \frac{1}{4}$, $i=1,2$. As previously, we generate, using $%
F_{i}$, the linear and continuous operators $R_{i}\in \mathcal{L}^{2}(X)$,
hence $\left\Vert R_{i}\right\Vert _{o}\leq Q_{i}\leq \frac{1}{4}$, $i=1,2$.
Then $\left\Vert R_{1}\right\Vert _{o}+\left\Vert R_{2}\right\Vert _{o}\leq 
\frac{1}{2}<1$.

Take also $\mu ^{0}\in cabv(L^{2}(\lambda ))$ with $\left\Vert \mu
^{0}\right\Vert \leq \frac{a}{2}$. Then 
\begin{equation*}
\left\Vert \mu ^{0}\right\Vert +a(\left\Vert R_{1}\right\Vert
_{o}+\left\Vert R_{2}\right\Vert _{o})\leq \frac{a}{2}+\frac{a}{2}=a\text{.}
\end{equation*}

Under these conditions, we can apply Theorem 2.2.3.

The effective computation will be done for the following particular case:

Take\ first $F_{1}(x,y)=\frac{xy}{2}$ and $F_{2}(x,y)=\frac{x^{2}y^{2}}{4}$,
hence $Q_{1}=Q_{2}=\frac{1}{4}$.

In order to introduce the measure $\mu ^{0}$, we consider first the measure $%
m\in cabv(L^{2}(\lambda ))$ given, for any $B\in \mathcal{B}$, via 
\begin{equation*}
m(B)=\overset{\sim }{h_{B}}\text{,}
\end{equation*}%
where $h_{B}:[0,1]\rightarrow K$ is the continuous function acting as
follows:%
\begin{equation*}
h_{B}(t)=\lambda (B\cap \lbrack 0,t])\text{,}
\end{equation*}%
$t\in \lbrack 0,1]$. Then, we know that $\left\Vert m\right\Vert =\frac{2}{3}
$ (see [5]). Finally, we take $\mu ^{0}\overset{def}{=}\frac{1}{2}m$, hence $%
\left\Vert \mu ^{0}\right\Vert =\frac{1}{3}$ and $a=1$. Consequently $%
\left\Vert R_{1}\right\Vert _{o}+\left\Vert R_{2}\right\Vert _{o}\leq
Q_{1}+Q_{2}=\frac{1}{2}<1$ and $\left\Vert \mu ^{0}\right\Vert +a(\left\Vert
R_{1}\right\Vert _{o}+\left\Vert R_{2}\right\Vert _{o})\leq \frac{1}{3}+%
\frac{1}{2}<1=a$.

We obtain the unique fractal (invariant) measure $\mu ^{\ast }\in
cabv(L^{2}(\lambda ))$.

The invariance equation is (for any $B\in \mathcal{B}$):

\begin{equation*}
R_{1}(\mu ^{\ast }((3B)\cap \lbrack 0,1]))+R_{2}(\mu ^{\ast }((3B-2)\cap
\lbrack 0,1]))+\mu ^{0}(B)=\mu ^{\ast }(B)\text{.}
\end{equation*}

Considering, for any $B\in \mathcal{B}$, a representative $\overset{\sim }{%
f_{B}}\in \mu ^{\ast }(B)$, we obtain representatives of $R_{1}(\mu ^{\ast
}(B))$, $R_{2}(\mu ^{\ast }(B))$ via 
\begin{equation*}
R_{1}(\mu ^{\ast }(B))=\frac{1}{4}\overset{1}{\underset{0}{\dint }}%
xyf_{B}(y)d\lambda (y)
\end{equation*}%
\begin{equation*}
R_{2}(\mu ^{\ast }(B))=\frac{1}{4}\overset{1}{\underset{0}{\dint }}%
x^{2}y^{2}f_{B}(y)d\lambda (y)
\end{equation*}%
and the invariance equation can be (abusively) written%
\begin{equation*}
\frac{1}{4}\overset{1}{\underset{0}{\dint }}xyf_{(3B)\cap \lbrack
0,1]}(y)d\lambda (y)+\frac{1}{4}\overset{1}{\underset{0}{\dint }}%
x^{2}y^{2}f_{(3B-2)\cap \lbrack 0,1]}(y)d\lambda (y)+\frac{1}{2}\lambda
(B\cap \lbrack 0,x])=f_{B}(x)
\end{equation*}%
for any $B\in \mathcal{B}$ and $\lambda $-almost all $x\in \lbrack 0,1]$.

In particular, for $B=[0,1]$, let us write $f_{[0,1]}=\varphi $, hence, for $%
\lambda $-almost all $x\in \lbrack 0,1]$, one has the integral equation%
\begin{equation*}
\varphi (x)=\frac{1}{2}x+\frac{1}{4}(x\overset{1}{\underset{0}{\dint }}%
y\varphi (y)d\lambda (y)+x^{2}\overset{1}{\underset{0}{\dint }}y^{2}\varphi
(y)d\lambda (y))
\end{equation*}%
and this gives%
\begin{equation*}
\varphi (x)=\frac{1}{2}x+\alpha x+\beta x^{2}\text{,}
\end{equation*}%
where $\alpha =\frac{1}{4}\overset{1}{\underset{0}{\dint }}y\varphi
(y)d\lambda (y)$ and $\beta =\frac{1}{4}\overset{1}{\underset{0}{\dint }}%
y^{2}\varphi (y)d\lambda (y))$.

Substituting:%
\begin{equation*}
\varphi (x)=\frac{1}{2}x+\alpha x+x^{2}=
\end{equation*}%
\begin{equation*}
=\frac{1}{2}x+\frac{1}{4}[x\overset{1}{\underset{0}{\dint }}y(\frac{1}{2}%
y+\alpha y+\beta y^{2})d\lambda (y)+x^{2}\overset{1}{\underset{0}{\dint }}%
y^{2}(\frac{1}{2}y+\alpha y+\beta y^{2})d\lambda (y)]=
\end{equation*}%
\begin{equation*}
=\frac{1}{2}x+\frac{1}{4}[x(\frac{1}{6}+\frac{\alpha }{3}+\frac{\beta }{4}%
)+x^{2}(\frac{1}{8}+\frac{\alpha }{4}+\frac{\beta }{5})]\text{.}
\end{equation*}

Identifying:%
\begin{equation*}
\{%
\begin{array}{c}
\alpha =\frac{1}{24}+\frac{\alpha }{12}+\frac{\beta }{16} \\ 
\beta =\frac{1}{32}+\frac{\alpha }{16}+\frac{B}{20}%
\end{array}%
\Leftrightarrow \{%
\begin{array}{c}
\frac{11}{12}\alpha -\frac{1}{16}\beta =\frac{1}{24} \\ 
-\frac{1}{16}\alpha +\frac{19}{20}\beta =\frac{1}{32}%
\end{array}%
\end{equation*}%
with solutions $\alpha =\frac{319}{6658}$ and $\beta =\frac{120}{3329}$.

Finally%
\begin{equation*}
\varphi (x)=\frac{24}{3329}(76x+5x^{2})\text{.}
\end{equation*}

\bigskip

2.3. \textbf{The Particular Case When All \ the Functions }$\omega _{i}$ 
\textbf{Are Constant}

\bigskip

We shall consider the initial framework of the paragraph, adding the
following supplementary:

\textit{Assumption}. All the functions $\omega _{i}$ are constant.

More precisely, we consider $M$ distinct points $t_{1}$, $t_{2}$, ..., $%
t_{M} $ in $T$ such that for any $i=1,2,...,M$ and any $t\in T$ one has $%
\omega _{i}(t)=t_{i}$.

Under this assumption, it is easy to see that, for any $i=1,2,...,M$ and any 
$\mu \in cabv(X)$, one has%
\begin{equation*}
\omega _{i}(\mu )(B)=\delta _{t_{i}}(B)\mu (T)\text{,}
\end{equation*}%
if $B\in \mathcal{B}$. This leads to the formula, valid for any $\mu \in
cabv(X)$:%
\begin{equation*}
H(\mu )=\overset{M}{\underset{i=1}{\sum }}\delta _{t_{i}}R_{i}(\mu (T))
\end{equation*}%
(here $\delta _{t_{i}}:\mathcal{B\rightarrow }\mathbb{R}_{+}$ is the Dirac
measure concentrated at $t_{i}$).

Define the set $T_{1}=\{t_{1},t_{2},...,t_{M}\}$. It follows that, for any $%
\mu \in cabv(X)$ and any $B\in \mathcal{B}$, one has%
\begin{equation*}
H(\mu )(B)=\{%
\begin{array}{cc}
\underset{t_{i}\in B}{\dsum }R_{i}(\mu (T))\text{,} & \text{if }B\cap
T_{1}\neq \emptyset \\ 
0\text{,} & \text{if }B\cap T_{1}=\emptyset%
\end{array}%
\end{equation*}%
which will be written in the sequel in all cases%
\begin{equation*}
H(\mu )(B)=\underset{t_{i}\in B}{\dsum }R_{i}(\mu (T))\text{,}
\end{equation*}%
hence%
\begin{equation*}
H(\mu )=\overset{M}{\underset{i=1}{\sum }}\delta _{t_{i}}R_{i}(\mu (T))
\end{equation*}%
and, for any $B\in \mathcal{B}$, 
\begin{equation*}
H(\mu )(B)=H(\mu )(B\cap T_{1})\text{.}
\end{equation*}

The reader can easily adapt theorems 2.2.1 and 2.2.2 to this particular case
(we have all $r_{i}=0$). We shall concentrate our attention to Theorem
2.2.3, because the particular form of this theorem in case all $\omega _{i}$
are constant is more interesting.

In order to proceed further, we recall that $(\mathcal{L}(X),\left\Vert
.\right\Vert _{o})$ is a Banach algebra with multiplication given by $UV%
\overset{def}{=}U\circ V$. For more details, see [1].

\bigskip

\textbf{Theorem 2.3.1.} \textit{Let us accept the present condition, i.e.
all }$\omega _{i}$\textit{\ are constant. Assume also that} $\mu ^{0}\in
cabv(X)$ \textit{and} 
\begin{equation*}
e=\overset{M}{\underset{i=1}{\sum }}\left\Vert R_{i}\right\Vert _{o}<1\text{.%
}
\end{equation*}

\textit{Define} $H_{2}:cabv(X)\rightarrow cabv(X)$\textit{\ via }%
\begin{equation*}
H_{2}(\mu )=H(\mu )+\mu ^{0}=\overset{M}{\underset{i=1}{\sum }}\delta
_{t_{i}}R_{i}(\mu (T))+\mu ^{0}\text{.}
\end{equation*}

\textit{It follows that} $H_{2}$\textit{\ is a contraction with contraction
factor }$\leq e$\textit{\ and possesses an unique fractal (invariant)
measure }$\mu ^{\ast }\in cabv(X)$\textit{, i.e. }$H_{2}(\mu ^{\ast })=\mu
^{\ast }$\textit{.}

\textit{Moreover, writing }$R\overset{def}{=}\overset{M}{\underset{i=1}{\sum 
}}R_{i}\in \mathcal{L}(X)$\textit{, it follows that }$Id_{X}-R$\textit{\ is
invertible in} $\mathcal{L}(X)$ \textit{and we have the formula}%
\begin{equation*}
\mu ^{\ast }=\overset{M}{\underset{i=1}{\sum }}\delta _{t_{i}}(R_{i}\circ
(Id_{X}-R)^{-1})(\mu ^{0}(T))+\mu ^{0}\text{.}
\end{equation*}

\textit{Proof}. The action of $H_{2}$ is correctly defined, in view of the
beginning of the paragraph and of the previous computations. The fact that $%
Id_{X}-R$\textit{\ }is invertible follows from $\left\Vert R\right\Vert
_{o}\leq \overset{M}{\underset{i=1}{\sum }}\left\Vert R_{i}\right\Vert
_{o}<1 $.

The existence and uniqueness of $\mu ^{\ast }$ follow from Theorem 2.2.3.

Notice that, because $H_{2}(\mu ^{\ast })=\mu ^{\ast }$, we have, for any $%
B\in \mathcal{B}$%
\begin{equation}
\mu ^{\ast }(B)=\mu ^{0}(B)+\underset{t_{i}\in B}{\dsum }R_{i}(\mu ^{\ast
}(T))\text{,}  \tag{2.3.1}
\end{equation}%
hence, for $B\cap T_{1}=\emptyset $ one has $\mu ^{\ast }(B)=\mu ^{0}(B)$
and, in particular%
\begin{equation}
\mu ^{\ast }(T\smallsetminus T_{1})=\mu ^{0}(T)-\mu ^{0}(T_{1})\text{.} 
\tag{2.3.2}
\end{equation}

Because, for any $i$, one has%
\begin{equation*}
\mu ^{\ast }(\{t_{i}\})=\mu ^{0}(\{t_{i}\})+R_{i}(\mu ^{\ast }(T))\text{,}
\end{equation*}%
it follows, by addition, that%
\begin{equation*}
\mu ^{\ast }(T_{1})=\mu ^{0}(T_{1})+(\overset{M}{\underset{i=1}{\sum }}%
R_{i})(\mu ^{\ast }(T_{1})+\mu ^{\ast }(T\smallsetminus T_{1}))
\end{equation*}%
which means (see 2.3.2) that%
\begin{equation*}
\mu ^{\ast }(T_{1})=\mu ^{0}(T_{1})+R(\mu ^{\ast }(T_{1})+\mu ^{0}(T)-\mu
^{0}(T_{1}))\Leftrightarrow
\end{equation*}%
\begin{equation*}
\Leftrightarrow (Id_{X}-R)(\mu ^{\ast }(T_{1}))=(Id_{X}-R)(\mu
^{0}(T_{1}))+R(\mu ^{0}(T))\Leftrightarrow
\end{equation*}%
\begin{equation*}
\Leftrightarrow (Id_{X}-R)(\mu ^{\ast }(T_{1})-\mu ^{0}(T_{1}))=R(\mu
^{0}(T))\Leftrightarrow
\end{equation*}%
\begin{equation*}
\Leftrightarrow \mu ^{\ast }(T_{1})-\mu ^{0}(T_{1})=((Id_{X}-R)^{-1}\circ
R)(\mu ^{0}(T))\text{.}
\end{equation*}

It follows (see (2.3.2)) that%
\begin{equation*}
\mu ^{\ast }(T)=\mu ^{\ast }(T_{1})+\mu ^{\ast }(T\smallsetminus T_{1})=\mu
^{\ast }(T_{1})+\mu ^{0}(T)\smallsetminus \mu ^{0}(T_{1})=
\end{equation*}%
\begin{equation*}
=\mu ^{0}(T)+(Id_{X}-R)^{-1}\circ R)(\mu ^{0}(T))=
\end{equation*}%
\begin{equation*}
=(Id_{X}+(Id_{X}-R)^{-1}\circ R)(\mu ^{0}(T))=(Id_{X}-R)^{-1}\circ R)(\mu
^{0}(T))\text{.}
\end{equation*}

Using (2.3.1), we finally get, for any $B\in \mathcal{B}$:%
\begin{equation}
\mu ^{\ast }(B)=\mu ^{0}(B)+\underset{t_{i}\in B}{\dsum }(R_{i}\circ
(Id_{X}-R)^{-1})(\mu ^{0}(T))  \tag{2.3.3}
\end{equation}%
and (2.3.3) proves the enunciation. $\square $

\bigskip

\textbf{Remark}. For any $\mu $ and $\mu ^{0}$ in $cabv(X)$, the measure 
\begin{equation*}
\mu ^{0}+\overset{M}{\underset{i=1}{\sum }}\delta _{t_{i}}R_{i}(\mu (T))
\end{equation*}%
is clearly in $cabv(X)$. Hence, one can define $H_{2}:cabv(X)\rightarrow
cabv(X)$ like in the enunciation of Theorem 2.3.1, without assuming that $%
e<1 $.

Moreover, if $\left\Vert R\right\Vert _{o}<1$ (it is sufficient to have $%
Id_{X}-R$ invertible), the formula defining $\mu ^{\ast }$ in the
enunciation of Theorem 2.3.1 gives a fixed point $\mu ^{\ast }\in cabv(X)$
of $H_{2}$: $H_{2}(\mu ^{\ast })=\mu ^{\ast }$ (easy computation, because $%
\mu ^{\ast }(T)=(Id_{X}-R)^{-1}(\mu ^{0}(T))$).

\bigskip

The following example will involve all the conditions in Theorem 2.3.1.

\bigskip

\textbf{Concrete Illustrations of Theorem 2.3.1}

\bigskip

We shall work for $T=[0,1]$, $X=K^{2}$, $M=2$ and $t_{1}=0$, $t_{2}=1$.

Take $\mu ^{0}\in cabv(K^{2})$, given via $\mu ^{0}(B)=(\lambda (B),\delta
_{0}(B))$ for any $B\in \mathcal{B}$ (where $\lambda :\mathcal{B}\rightarrow 
\mathbb{R}_{+}$ is the Lebesgue measure on $[0,1]$ and $\delta _{0}:\mathcal{%
B}\rightarrow \mathbb{R}_{+}$ is the Dirac measure concentrated at $0$).

For $R_{1}$,$R_{1}\in \mathcal{L}(K^{2})$ given via%
\begin{equation*}
R_{1}\equiv (%
\begin{array}{cc}
\frac{1}{8} & \frac{1}{8} \\ 
-\frac{1}{16} & \frac{1}{8}%
\end{array}%
)\text{, }R_{2}\equiv (%
\begin{array}{cc}
\frac{1}{8} & -\frac{1}{8} \\ 
\frac{1}{16} & \frac{1}{8}%
\end{array}%
)
\end{equation*}%
one has 
\begin{equation*}
\left\Vert R_{1}\right\Vert _{o}=\left\Vert R_{2}\right\Vert _{o}=\frac{3}{16%
}\text{ with }\left\Vert R_{1}\right\Vert _{o}+\left\Vert R_{2}\right\Vert
_{o}=\frac{3}{8}<1
\end{equation*}%
and 
\begin{equation*}
R=R_{1}+R_{2}\equiv (%
\begin{array}{cc}
\frac{1}{4} & 0 \\ 
0 & \frac{1}{4}%
\end{array}%
)\Rightarrow Id_{K^{2}}-R\equiv (%
\begin{array}{cc}
\frac{3}{4} & 0 \\ 
0 & \frac{3}{4}%
\end{array}%
)
\end{equation*}%
with%
\begin{equation*}
(Id_{K^{2}}-R)^{-1}\equiv (%
\begin{array}{cc}
\frac{4}{3} & 0 \\ 
0 & \frac{4}{3}%
\end{array}%
)\text{.}
\end{equation*}

We have $\mu ^{0}(T)=(1,1)$ and%
\begin{equation*}
R_{1}\circ (Id_{K^{2}}-R)^{-1}\equiv (%
\begin{array}{cc}
\frac{1}{8} & \frac{1}{8} \\ 
-\frac{1}{16} & \frac{1}{8}%
\end{array}%
)(%
\begin{array}{cc}
\frac{4}{3} & 0 \\ 
0 & \frac{4}{3}%
\end{array}%
)=(%
\begin{array}{cc}
\frac{1}{6} & \frac{1}{6} \\ 
-\frac{1}{12} & \frac{1}{6}%
\end{array}%
)
\end{equation*}%
\begin{equation*}
R_{2}\circ (Id_{K^{2}}-R)^{-1}\equiv (%
\begin{array}{cc}
\frac{1}{8} & -\frac{1}{8} \\ 
\frac{1}{16} & \frac{1}{8}%
\end{array}%
)(%
\begin{array}{cc}
\frac{4}{3} & 0 \\ 
0 & \frac{4}{3}%
\end{array}%
)=(%
\begin{array}{cc}
\frac{1}{6} & -\frac{1}{6} \\ 
\frac{1}{12} & \frac{1}{6}%
\end{array}%
)\text{.}
\end{equation*}

The theorem works. Take $B\in \mathcal{B}$ and let us compute $\mu ^{\ast
}(B)$ (see (2.3.3)):

a) If $0\in B$ and $1\in B$:%
\begin{equation*}
\mu ^{\ast }(B)\equiv (%
\begin{array}{c}
\lambda (B) \\ 
1%
\end{array}%
)+(%
\begin{array}{cc}
\frac{1}{6} & \frac{1}{6} \\ 
-\frac{1}{12} & \frac{1}{6}%
\end{array}%
)(%
\begin{array}{c}
1 \\ 
1%
\end{array}%
)+(%
\begin{array}{cc}
\frac{1}{6} & -\frac{1}{6} \\ 
\frac{1}{12} & \frac{1}{6}%
\end{array}%
)(%
\begin{array}{c}
1 \\ 
1%
\end{array}%
)=
\end{equation*}%
\begin{equation*}
=(%
\begin{array}{c}
\lambda (B) \\ 
1%
\end{array}%
)+(%
\begin{array}{c}
\frac{1}{3} \\ 
\frac{1}{3}%
\end{array}%
)=(%
\begin{array}{c}
\lambda (B)+\frac{1}{3} \\ 
\frac{4}{3}%
\end{array}%
)\text{.}
\end{equation*}

b) If $0\in B$ and $1\notin B$:%
\begin{equation*}
\mu ^{\ast }(B)\equiv (%
\begin{array}{c}
\lambda (B) \\ 
1%
\end{array}%
)+(%
\begin{array}{cc}
\frac{1}{6} & \frac{1}{6} \\ 
-\frac{1}{12} & \frac{1}{6}%
\end{array}%
)(%
\begin{array}{c}
1 \\ 
1%
\end{array}%
)=
\end{equation*}%
\begin{equation*}
=(%
\begin{array}{c}
\lambda (B) \\ 
1%
\end{array}%
)+(%
\begin{array}{c}
\frac{1}{3} \\ 
\frac{1}{12}%
\end{array}%
)=(%
\begin{array}{c}
\lambda (B)+\frac{1}{3} \\ 
\frac{13}{12}%
\end{array}%
)\text{.}
\end{equation*}

c) If $0\notin B$ and $1\in B$:%
\begin{equation*}
\mu ^{\ast }(B)\equiv (%
\begin{array}{c}
\lambda (B) \\ 
0%
\end{array}%
)+(%
\begin{array}{cc}
\frac{1}{6} & -\frac{1}{6} \\ 
\frac{1}{12} & \frac{1}{6}%
\end{array}%
)(%
\begin{array}{c}
1 \\ 
1%
\end{array}%
)=
\end{equation*}%
\begin{equation*}
=(%
\begin{array}{c}
\lambda (B) \\ 
0%
\end{array}%
)+(%
\begin{array}{c}
0 \\ 
\frac{1}{4}%
\end{array}%
)=(%
\begin{array}{c}
\lambda (B) \\ 
\frac{1}{4}%
\end{array}%
)\text{.}
\end{equation*}

d) If $0\notin B$ and $1\notin B$:%
\begin{equation*}
\mu ^{\ast }(B)\equiv (%
\begin{array}{c}
\lambda (B) \\ 
0%
\end{array}%
)\text{.}
\end{equation*}

\bigskip

\textbf{3. Fractal (Invariant) Vector Measures. The Countable Case}

\bigskip

The idea is to replace the finite sets $\{R_{i}\mid i=1,2,...,M\}$ and $%
\{\omega _{i}\mid i=1,2,...,M\}$ with countable sets $\{R_{i}\mid i\in 
\mathbb{N}^{\ast }\}$ and $\{\omega _{i}\mid i\in \mathbb{N}^{\ast }\}$ and
to obtain similar results. Proofs will be, many times, skipped or merely
sketched, laying stress upon the facts in the proofs which differ
essentially from those in the finite case.

\bigskip

3.1 \textbf{Framework of the Paragraph}

\bigskip

Again $(T,d)$ is a compact metric space with Borel sets $\mathcal{B}$ and $X$
is a Hilbert space.

We consider a generalized iterated function system, i.e. a sequence $(\omega
_{i})_{i\geq 1}\subset Lip(T)$ with Lipschitz constants $r_{i}=\left\Vert
\omega _{i}\right\Vert _{L}$ and a sequence $(R_{i})_{i\geq 1}\subset 
\mathcal{L}(X)$.

We assume that the sequence $(r_{i})_{i\geq 1}$ is bounded and the series $%
\overset{\infty }{\underset{i=1}{\dsum }}\left\Vert R_{i}\right\Vert _{o}$
is convergent.

We notice that, for any $\mu \in cabv(X)$, the series $\overset{\infty }{%
\underset{i=1}{\dsum }}R_{i}\circ \omega _{i}(\mu )$ is absolutely
convergent in $cabv(X)$.

Indeed, for any $i$ one has (norm in $cabv(X))$):%
\begin{equation*}
\left\Vert R_{i}\circ \omega _{i}(\mu )\right\Vert =\left\vert R_{i}\circ
\omega _{i}(\mu )\right\vert (T)\leq \left\Vert R_{i}\right\Vert
_{o}\left\vert \omega _{i}(\mu )\right\vert (T)=\left\Vert R_{i}\right\Vert
_{o}\left\Vert \omega _{i}(\mu )\right\Vert \leq \left\Vert R_{i}\right\Vert
_{o}\left\Vert \mu \right\Vert \text{.}
\end{equation*}

Denote by $H(\mu )\in cabv(X)$ the sum of this series, hence we have, for
any $\mu \in cabv(X)$%
\begin{equation}
H(\mu )=\overset{\infty }{\underset{i=1}{\dsum }}R_{i}\circ \omega _{i}(\mu )%
\text{.}  \tag{3.1.1}
\end{equation}

Of course, this relation is valid pointwise too, hence, for any $B\in 
\mathcal{B}$%
\begin{equation*}
H(\mu )(B)=\overset{\infty }{\underset{i=1}{\dsum }}(R_{i}\circ \omega
_{i}(\mu ))(B)=\overset{\infty }{\underset{i=1}{\dsum }}R_{i}(\mu (\omega
_{i}^{-1}(B)))\text{,}
\end{equation*}%
the last series being absolutely convergent in $X$. Clearly 
\begin{equation*}
\left\Vert H(\mu )\right\Vert =\underset{M}{\lim }\left\Vert \overset{M}{%
\underset{i=1}{\dsum }}R_{i}\circ \omega _{i}(\mu )\right\Vert
\end{equation*}%
and, using the fact that $\left\Vert \overset{M}{\underset{i=1}{\dsum }}%
R_{i}\circ \omega _{i}(\mu )\right\Vert \leq (\overset{M}{\underset{i=1}{%
\dsum }}\left\Vert R_{i}\right\Vert _{o})\left\Vert \mu \right\Vert $ for
any $M$ (see the beginning of the preceding paragraph) we get%
\begin{equation*}
\left\Vert H(\mu )\right\Vert \leq (\overset{\infty }{\underset{i=1}{\dsum }}%
\left\Vert R_{i}\right\Vert _{o})\left\Vert \mu \right\Vert \text{.}
\end{equation*}

We proved that formula (3.1.1) defines a linear and continuous operator $%
H\in \mathcal{L}(cabv(X))$ with%
\begin{equation}
\left\Vert H\right\Vert _{o}\leq \overset{\infty }{\underset{i=1}{\dsum }}%
\left\Vert R_{i}\right\Vert _{o}\text{.}  \tag{3.1.2}
\end{equation}

Like in the finite case, we call $H$ the Markov-type operator (generated by $%
(R_{i})_{i\geq 1}$ and $(\omega _{i})_{i\geq 1}$).

In the sequel, we follow the same lines as in the finite case.

\bigskip

\textbf{Lemma 3.1.1. }\textit{Let\ }$f\in L_{1}(T,X)$\textit{. Then:}

\textit{i) For any }$t\in T$\textit{, the series }$\overset{\infty }{%
\underset{i=1}{\sum }}(R_{i}^{\ast }\circ f\circ \omega _{i})(t)$\textit{\
is absolutely convergent.}

\textit{ii) Define the function }$g:T\rightarrow X$\textit{, via }$g(t)=%
\overset{\infty }{\underset{i=1}{\sum }}(R_{i}^{\ast }\circ f\circ \omega
_{i})(t)$\textit{\ for any }$t\in T$\textit{\ (according to point (i)). Then 
}$g\in Lip(T,X)$\textit{\ and}%
\begin{equation*}
\left\Vert g\right\Vert _{L}\leq \overset{\infty }{\underset{i=1}{\sum }}%
\left\Vert R_{i}\right\Vert _{o}r_{i}<\infty \text{.}
\end{equation*}

\textit{Proof}. For any $t\in T$, one has%
\begin{equation*}
\overset{\infty }{\underset{i=1}{\sum }}\left\Vert R_{i}^{\ast }(f(\omega
_{i}((t)))\right\Vert \leq \overset{\infty }{\underset{i=1}{\sum }}%
\left\Vert R_{i}^{\ast }\right\Vert _{o}\left\Vert f(\omega
_{i}((t))\right\Vert \leq \overset{\infty }{\underset{i=1}{\sum }}\left\Vert
R_{i}\right\Vert _{o}\left\Vert f\right\Vert _{\infty }<\infty \text{,}
\end{equation*}%
thus proving (i).

For (ii):%
\begin{equation*}
\left\Vert g(x)-g(y)\right\Vert \leq \overset{\infty }{\underset{i=1}{\sum }}%
\left\Vert R_{i}^{\ast }\right\Vert _{o}\left\Vert f(\omega
_{i}(x))-f(\omega _{i}(y))\right\Vert \leq (\overset{M}{\underset{i=1}{\sum }%
}\left\Vert R_{i}\right\Vert _{o}r_{i})d(x,y)\text{. }\square
\end{equation*}

\bigskip

\textbf{Theorem 3.1.2}. \textit{(Change of Variable Formula). For any }$f\in
C(X)$\textit{\ and any }$\mu \in cabv(X)$\textit{, one has}%
\begin{equation*}
\dint fdH(\mu )=\dint gd\mu \text{,}
\end{equation*}%
\textit{where }$g=\overset{\infty }{\underset{i=1}{\sum }}R_{i}^{\ast }\circ
f\circ \omega _{i}$\textit{.}

\textit{Proof}. For any $M\in \mathbb{N}^{\ast }$, we have, according to
Theorem 2.1.2:%
\begin{equation*}
\dint fdH_{M}(\mu )=\dint g_{M}d\mu \text{,}
\end{equation*}%
where $H_{M}(\mu )\overset{def}{=}\overset{M}{\underset{i=1}{\sum }}%
R_{i}\circ \omega _{i}(\mu )$ and $g_{M}\overset{def}{=}\overset{M}{\underset%
{i=1}{\sum }}R_{i}^{\ast }\circ f\circ \omega _{i}$.

Using (3.1.3), it will be sufficient to prove that%
\begin{equation*}
\underset{M}{\lim }\dint fdH_{M}(\mu )=\dint fdH(\mu )\text{ and }\underset{M%
}{\lim }\dint g_{M}d\mu =\dint gd\mu \text{.}
\end{equation*}

Because of the fact (valid for any $h\in C(X)$ and any $m\in cabv(X)$) that $%
\left\vert \int hdm\right\vert \leq \left\Vert h\right\Vert _{\infty
}\left\Vert m\right\Vert $, it will be sufficient to prove that $H_{M}(\mu )%
\underset{M}{\rightarrow }H(\mu )$ in $cabv(X)$ and $g_{M}\underset{M}{%
\overset{u}{\rightarrow }}g$.

To this end, fix an arbitrary $\varepsilon >0$ and find $M_{0}\in \mathbb{N}%
^{\ast }$ such that $\overset{\infty }{\underset{i=M+1}{\sum }}\left\Vert
R_{i}\right\Vert _{o}<\varepsilon $, whenever $M\geq M_{0}$.

First take $M\geq M_{0}$ and a partition $(A_{j})_{1\leq j\leq n}$. Consider
the measure $P_{M}=H(\mu )-H_{M}(\mu )$ which acts pointwise via $P_{M}(B)=%
\overset{\infty }{\underset{i=M+1}{\sum }}R_{i}(\mu (\omega _{i}^{-1}(B)))$.
Then%
\begin{equation*}
\overset{n}{\underset{j=1}{\sum }}\left\Vert P(A_{j})\right\Vert \leq 
\overset{n}{\underset{j=1}{\sum }}\overset{\infty }{\underset{i=M+1}{\sum }}%
\left\Vert R_{i}(\mu (\omega _{i}^{-1}(A_{j})))\right\Vert \leq \overset{n}{%
\underset{j=1}{\sum }}\overset{\infty }{\underset{i=M+1}{\sum }}\left\Vert
R_{i}\right\Vert _{o}\left\Vert \mu (\omega _{i}^{-1}(A_{j}))\right\Vert =
\end{equation*}%
\begin{equation*}
=\overset{\infty }{\underset{i=M+1}{\sum }}\left\Vert R_{i}\right\Vert _{o}%
\overset{n}{\underset{j=1}{\sum }}\left\Vert \mu (\omega
_{i}^{-1}(A_{j}))\right\Vert \leq \overset{\infty }{\underset{i=M+1}{\sum }}%
\left\Vert R_{i}\right\Vert _{o}\left\Vert \mu \right\Vert \leq \varepsilon
\left\Vert \mu \right\Vert \text{.}
\end{equation*}%
Hence $P_{M}\underset{M}{\rightarrow }0$ in $cabv(X)$, i.e. $H_{M}(\mu )%
\underset{M}{\rightarrow }H(\mu )$.

Next take $M\geq M_{0}$ and $t\in T$: 
\begin{equation*}
\left\Vert g(t)-g_{M}(t)\right\Vert =\left\Vert \overset{\infty }{\underset{%
i=M+1}{\sum }}R_{i}^{\ast }(f(\omega _{i}(t)))\right\Vert \leq \overset{%
\infty }{\underset{i=M+1}{\sum }}\left\Vert R_{i}^{\ast }\right\Vert
_{o}\left\Vert f(\omega _{i}(t))\right\Vert \leq
\end{equation*}%
\begin{equation*}
\leq \overset{\infty }{\underset{i=M+1}{\sum }}\left\Vert R_{i}^{\ast
}\right\Vert _{o}\left\Vert f\right\Vert _{\infty }\leq \varepsilon
\left\Vert f\right\Vert _{\infty }\text{,}
\end{equation*}%
i.e. $g_{M}\underset{M}{\overset{u}{\rightarrow }}g$. $\square $

\bigskip

\textbf{Theorem 3.1.3.} \textit{Considering on }$cabv(X)$\textit{\ the norm }%
$\left\Vert .\right\Vert _{MK}$\textit{, we have }$H\in \mathcal{L}(cabv(X))$%
\textit{\ and}%
\begin{equation*}
\left\Vert H\right\Vert _{o}\leq \overset{\infty }{\underset{i=1}{\dsum }}%
\left\Vert R_{i}\right\Vert _{o}(1+r_{i})<\infty \text{.}
\end{equation*}

\textit{Sketch of proof}. We follow the lines in the proof of Theorem 2.1.3.

Namely, we construct $g$ as in Lemma 3.1.1, and notice that $\left\Vert
g\right\Vert _{\infty }\leq \overset{\infty }{\underset{i=1}{\dsum }}%
\left\Vert R_{i}\right\Vert _{o}$, $\left\Vert g\right\Vert _{L}\leq \overset%
{\infty }{\underset{i=1}{\dsum }}\left\Vert R_{i}\right\Vert _{o}r_{i}$,
hence $\left\Vert g\right\Vert _{BL}\leq \overset{\infty }{\underset{i=1}{%
\dsum }}\left\Vert R_{i}\right\Vert _{o}(1+r_{i})$. Finally, we use Theorem
3.1.2. $\square $

\bigskip

\textbf{Theorem 3.1.4.} \textit{For any} $\mu \in cabv(X,0)$\textit{, one
has }$H(\mu )\in cabv(X,0)$\textit{. Hence, one can define }$%
H_{0}:cabv(X,0)\rightarrow cabv(X,0)$\textit{, via }$H_{0}(\mu )=H(\mu )$%
\textit{.}

\textit{Considering on }$cabv(X,0)$\textit{\ the norm }$\left\Vert
.\right\Vert _{MK}^{\ast }$\textit{, we have\ }$H_{0}\in \mathcal{L}%
(cabv(X,0))$\textit{\ and }%
\begin{equation*}
\left\Vert H_{0}\right\Vert _{o}\leq \overset{M}{\underset{i=1}{\dsum }}%
\left\Vert R_{i}\right\Vert _{o}r_{i}\text{.}
\end{equation*}

\textit{Sketch of proof}. We follow the lines in the proof of Theorem 2.1.4.

We construct $g$ as in Lemma 3.1.1, and notice that $\left\Vert g\right\Vert
_{L}\leq \overset{\infty }{\underset{i=1}{\dsum }}\left\Vert
R_{i}\right\Vert _{o}r_{i}$. Finally, we use Theorem 3.1.2. $\square $

\bigskip

The results in this subparagraph will be used in the next subparagraph,
where we shall study (like in the finite case) the $H_{1}$ and $H_{2}$
models.

\newpage

3.2. \textbf{Fixed Point Models}

\bigskip

First, we consider the operator $R\in \mathcal{L}(X)$, defined pointwise via%
\begin{equation*}
R=\overset{\infty }{\underset{i=1}{\dsum }}R_{i}
\end{equation*}%
(the convergence in $\mathcal{L}(X)$ is absolute).

The next three theorems are proved exactly like theorems 2.2.1, 2.2.2 and
2.2.3 (with minor adaptations in the proofs).

\bigskip

\textbf{Theorem 3.2.1.} \textit{Consider }$X=K^{n}$\textit{, }$n\in \mathbb{N%
}^{\ast }$\textit{. The hypotheses are:}

\textit{a) }$R=Id_{K^{n}}$;

\textit{b) }$c=\overset{\infty }{\underset{i=1}{\sum }}\left\Vert
R_{i}\right\Vert _{o}r_{i}<1$\textit{\ (Clearly this true if all }$\omega
_{i}$\textit{\ are contractions and }$\overset{\infty }{\underset{i=1}{\sum }%
}\left\Vert R_{i}\right\Vert _{o}=1$.)\textit{;}

\textit{c) }$0<a<\infty $\textit{\ and }$\upsilon \in K^{n}$\textit{\ are
such that }$\left\Vert \upsilon \right\Vert \leq a$ \textit{(hence }$%
B_{a}(K^{n},v)\neq \emptyset $\textit{);}

\textit{d) }$\emptyset \neq A\subseteq B_{a}(K^{n},v)$\textit{\ is such that 
}$H(A)\subseteq A$\textit{\ and }$A$\textit{\ is weak}$^{\ast }$\textit{\
closed (In the particular case when }$\left\Vert H(\mu )\right\Vert \leq
\left\Vert \mu \right\Vert $\textit{\ for any }$\mu \in cabv(K^{n})$\textit{%
, one can take }$A=B_{a}(K^{n},v)$\textit{. More particular, if }$\overset{%
\infty }{\underset{i=1}{\sum }}\left\Vert R_{i}\right\Vert _{o}=1$\textit{,
it follows that }$\left\Vert H(\mu )\right\Vert \leq \left\Vert \mu
\right\Vert $\textit{\ for any }$\mu \in cabv(K^{n})$.\textit{).}

\textit{Under these hypotheses, we define }$H_{1}:A\rightarrow A$\textit{\
via }$H_{1}(\mu )=H(\mu )$\textit{\ for any }$\mu \in A$\textit{. It follows
that }$H_{1}$\textit{\ is a contraction with contraction factor }$\leq c$%
\textit{, if }$A$\textit{\ is equipped with the metric }$d_{MK}^{\ast }$%
\textit{\ given via }$d_{MK}^{\ast }(\mu ,\nu )=\left\Vert \mu -\nu
\right\Vert _{MK}^{\ast }$\textit{.}

\textit{There exists an unique fractal (invariant) measure }$\mu ^{\ast }\in
A$\textit{\ of }$H_{1}$\textit{, i.e. }$H_{1}(\mu ^{\ast })=\mu ^{\ast }$%
\textit{.}

\bigskip

\textbf{Remarks}

1. Again, because $R=Id_{K^{n}}$, we have $1=\left\Vert R\right\Vert
_{o}\leq \overset{\infty }{\underset{i=1}{\sum }}\left\Vert R_{i}\right\Vert
_{o}$ and the condition $\overset{\infty }{\underset{i=1}{\sum }}\left\Vert
R_{i}\right\Vert _{o}=1$ is extremal.

2. Taking in the "classical model" a sequence $(p_{i})_{i\geq 1}$ with $%
p_{i}>0$ and $\overset{\infty }{\underset{i=1}{\sum }}p_{i}=1$ (instead $%
p_{1},p_{2},...,p_{M}$ with $\overset{M}{\underset{i=1}{\sum }}p_{i}=1$) we
obtain a situation when all the particular conditions are fulfilled.

\bigskip

\textbf{Theorem 3.2.2.} \textit{Consider }$X=K^{n}$\textit{, }$n\in \mathbb{N%
}^{\ast }$\textit{. The hypotheses are:}

\textit{a) }$d=\overset{\infty }{\underset{i=1}{\sum }}\left\Vert
R_{i}\right\Vert _{o}(1+r_{i})<1$\textit{\ (This true, in particular, if all 
}$\omega _{i}$\textit{\ are contractions and }$\overset{\infty }{\underset{%
i=1}{\sum }}\left\Vert R_{i}\right\Vert _{o}\leq \frac{1}{2})$\textit{;}

\textit{b) }$0<a<\infty $, $\mu ^{0}\in cabv(K^{n})$, $\emptyset \neq
A\subseteq B_{a}(K^{n})$\textit{\ is weak}$^{\ast }$\textit{\ closed and one
has }$H(\mu )+\mu ^{0}\in A$ \textit{for any} $\mu \in A$ \textit{(In
particular, if }$\left\Vert \mu ^{0}\right\Vert +a(\overset{\infty }{%
\underset{i=1}{\sum }}\left\Vert R_{i}\right\Vert _{o})\leq a$,\textit{\ one
can take }$A=B_{a}(K^{n})$.\textit{).}

\textit{Under these hypotheses, we define }$H_{2}:A\rightarrow A$\textit{\
via }$H_{2}(\mu )=H(\mu )+\mu ^{0}$\textit{\ for any }$\mu \in A$\textit{.
It follows that }$H_{2}$\textit{\ is a contraction with contraction factor }$%
\leq d$\textit{, if }$A$\textit{\ is equipped with the metric }$d_{MK}$%
\textit{, given via }$d_{MK}(\mu ,\nu )=\left\Vert \mu -\nu \right\Vert
_{MK} $\textit{.}

\textit{Then:}

\textit{i) If }$\mu ^{0}=0$\textit{, it follows that }$0\in A$\textit{.}

\textit{ii) There exists an unique fractal (invariant) measure }$\mu ^{\ast
}\in A$\textit{, i.e. }$H_{2}(\mu ^{\ast })=\mu ^{\ast }$\textit{. In case }$%
\mu ^{0}=0$\textit{, we have }$\mu ^{\ast }=0$.

\bigskip

\textbf{Theorem 3.2.3.} \textit{We work in an arbitrary Hilbert space }$X$ 
\textit{(as a matter of fact, this theorem is valid for any Banach space }$X$%
\textit{) and consider the usual Banach space }$cabv(X)$ \textit{with the
variational norm. Take }$\mu ^{0}\in cabv(X)$.

\textit{The hypotheses are:}

\textit{a) }$e=\overset{\infty }{\underset{i=1}{\sum }}\left\Vert
R_{i}\right\Vert _{o}<1$\textit{;}

\textit{b)} $\emptyset \neq A\subseteq cabv(X)$\textit{\ is a\ closed set
such that }$H(\mu )+\mu ^{0}\in A$ \textit{for any} $\mu \in A$ \textit{(In
particular, one can take }$A=cabv(X)$\textit{\ or one can take }$A=B_{a}(X)$%
, \textit{if }$0<a<\infty $ and $\left\Vert \mu ^{0}\right\Vert +a(\overset{%
\infty }{\underset{i=1}{\sum }}\left\Vert R_{i}\right\Vert _{o})\leq a$.%
\textit{).}

\textit{Under these hypotheses, define }$H_{2}:A\rightarrow A$\textit{\ via }%
$H_{2}(\mu )=H(\mu )+\mu ^{0}$\textit{\ for any }$\mu \in A$\textit{. It
follows that }$H_{2}$\textit{\ is a contraction with contraction factor }$%
\leq e$\textit{.}

\textit{Then:}

\textit{i) If }$\mu ^{0}=0$\textit{, then }$0\in A$\textit{.}

\textit{ii) There exists an unique fractal (invariant) measure }$\mu ^{\ast
}\in A$\textit{, i.e. }$H_{2}(\mu ^{\ast })=\mu ^{\ast }$\textit{. In case }$%
\mu ^{0}=0$\textit{, we have }$\mu ^{\ast }=0$.

\bigskip

We finish with the study of the particular case when all the contractions $%
\omega _{i}$ are constant, i.e. $\omega _{i}(t)=t_{i}\in T$ for any $t\in T$
and any $i\in \mathbb{N}^{\ast }$ (the points $t_{i}$ being distinct).

In this case, we have, for any $\mu \in cabv(X)$, the formula%
\begin{equation*}
H(\mu )=\overset{\infty }{\underset{i=1}{\sum }}\delta _{t_{i}}R_{i}(\mu (T))
\end{equation*}%
the convergence being in $cabv(X)$, due to the fact that $\omega _{i}(\mu
)=\delta _{t_{i}}\mu (T)$ for any $i$.

Again, theorems 3.2.1 and 3.2.2 can be easily adapted and we present the
adaptation of Theorem 3.2.3.

\bigskip

\textbf{Theorem 3.2.4.} \textit{Assume that all }$\omega _{i}$\textit{\ are
constant as previously. Let }$\mu ^{0}\in cabv(X)$\textit{\ and assume that}

\begin{equation*}
e=\overset{\infty }{\underset{i=1}{\sum }}\left\Vert R_{i}\right\Vert _{o}<1%
\text{\textit{.}}
\end{equation*}

\textit{Define} $H_{2}:cabv(X)\rightarrow cabv(X)$\textit{\ via }%
\begin{equation*}
H_{2}(\mu )=H(\mu )+\mu ^{0}=\overset{\infty }{\underset{i=1}{\sum }}\delta
_{t_{i}}R_{i}(\mu (T))+\mu ^{0}\text{.}
\end{equation*}

\textit{Then:}

\textit{i) }$Id_{X}-R$\textit{\ is invertible in }$\mathcal{L}(X)$.

\textit{ii)} $H_{2}$ \textit{is a contraction with contraction factor }$\leq
e$\textit{\ and possesses an unique fractal invariant (measure) }$\mu ^{\ast
}\in cabv(X)$\textit{, i.e. }$H_{2}(\mu ^{\ast })=\mu ^{\ast }$\textit{.
Moreover we have the formula}%
\begin{equation*}
\mu ^{\ast }=\overset{\infty }{\underset{i=1}{\sum }}\delta
_{t_{i}}(R_{i}\circ (Id_{X}-R)^{-1})(\mu ^{0}(T))+\mu ^{0}\text{.}
\end{equation*}

\textit{Sketch of proof}.

i) We have $\left\Vert R\right\Vert _{o}\leq \overset{\infty }{\underset{i=1}%
{\sum }}\left\Vert R_{i}\right\Vert _{o}<1$ and $Id_{X}-R$\textit{\ }is
invertible.

ii) We use Theorem 3.2.3. As concerns the formula for $\mu ^{\ast }$, we
notice that the formula in the present enunciation is meaningful, because,
for any $B\in \mathcal{B}$:%
\begin{equation*}
\overset{\infty }{\underset{i=1}{\sum }}\left\Vert \delta
_{t_{i}}(B)R_{i}((Id_{X}-R)^{-1}(\mu ^{0}(T)))\right\Vert \leq \overset{%
\infty }{\underset{i=1}{\sum }}\left\Vert R_{i}((Id_{X}-R)^{-1}(\mu
^{0}(T)))\right\Vert \leq
\end{equation*}%
\begin{equation*}
\leq (\overset{\infty }{\underset{i=1}{\sum }}\left\Vert R_{i}\right\Vert
_{o})\left\Vert (Id_{X}-R)^{-1}(\mu ^{0}(T))\right\Vert <\infty \text{.}
\end{equation*}

Then, we inspect the proof of\ Theorem 2.3.1 and notice that also in the
present case one has%
\begin{equation*}
\mu ^{\ast }(T_{1})-\mu ^{0}(T)=((Id_{X}-R)^{-1}\circ R)(\mu ^{0}(T))
\end{equation*}%
leading to%
\begin{equation*}
\mu ^{\ast }(T)=(Id_{X}-R)^{-1}(\mu ^{0}(T))
\end{equation*}%
and this finally proves that $H_{2}(\mu ^{\ast })=\mu ^{\ast }$. $\square $

\bigskip

\textbf{Remark. }Considerations similar to those in the Remark following
Theorem 2.3.1 are valid also here. So, one can construct $\mu ^{\ast }$
using the formula in the enunciation of Theorem 3.2.4 accepting (it is
sufficient to have $Id_{X}-R$ invertible) only the condition $\left\Vert
R\right\Vert _{o}<1$.

\bigskip

The following example illustrates Theorem 3.2.4 in the spirit of this Remark.

\bigskip

\textbf{Example 3.2.5}

\bigskip

Let $P\in \mathcal{L}(X)$ be arbitrary. In our schema (Theorem 3.2.4 plus
Remark), we shall take 
\begin{equation*}
R_{i}=-\frac{1}{i!}P^{i}
\end{equation*}%
$i=1,2,...$, working in the Banach algebra $\mathcal{L}(X)$. The
construction works, because 
\begin{equation*}
\overset{\infty }{\underset{i=1}{\sum }}\left\Vert R_{i}\right\Vert _{o}\leq 
\overset{\infty }{\underset{i=1}{\sum }}\frac{1}{i!}\left\Vert P\right\Vert
_{o}^{i}<\exp (\left\Vert P\right\Vert _{o})<\infty \text{.}
\end{equation*}

Defining $R=\overset{\infty }{\underset{i=1}{\sum }}R_{i}$, we see that 
\begin{equation*}
Id_{X}-R=Id_{X}+\overset{\infty }{\underset{i=1}{\sum }}\frac{1}{i!}%
P^{i}=\exp (P)\text{,}
\end{equation*}%
hence $Id_{X}-R$ is invertible with%
\begin{equation*}
(Id_{X}-R)^{-1}=\exp (-P)\text{.}
\end{equation*}

So, we have $H_{2}:cabv(X)\rightarrow cabv(X)$ given via%
\begin{equation*}
H_{2}(\mu )=-\overset{\infty }{\underset{i=1}{\sum }}\frac{1}{i!}\delta
_{t_{i}}P^{i}(\mu (T))+\mu ^{0}
\end{equation*}%
and $H_{2}$ possesses the fixed point $\mu ^{\ast }\in cabv(X)$ (i.e. $%
H_{2}(\mu ^{\ast })=\mu ^{\ast }$) given via%
\begin{equation*}
\mu ^{\ast }=-\overset{\infty }{\underset{i=1}{\sum }}\frac{1}{i!}\delta
_{t_{i}}(P^{i}\circ \exp (-P))(\mu ^{0}(T))+\mu ^{0}\text{.}
\end{equation*}

\bigskip

\textbf{References}

\bigskip

[1] G. R. Allan, Introduction to Banach Spaces and Algebras, Oxford
University Press, 2011.

[2] M. F. Barnsley, Fractals Everywhere (second edition), Morgan Kaufmann,
1993.

[3] K. Baron, A. Lasota, Markov operators on the space of vector measures;
coloured fractals, Ann. Pol. Math., 69 (1998), 217-234.

[4] I. Chi\c{t}escu, L. Ni\c{t}\u{a}, Fractal Vector Measures, Sci. Bull.,
Ser. A, Appl. Math. Phys., Politeh. Univ. Buchar., 77 (2015), 219-228.

[5] I. Chi\c{t}escu, L. Ioana, R. Miculescu, L. Ni\c{t}\u{a}, Sesquilinear
uniform integral, Proc. Indian Acad. Sci., Math. Sci., 125 (2015), 187-198.

[6] I. Chi\c{t}escu, L. Ioana, R. Miculescu, L. Ni\c{t}\u{a},
Monge-Kantorovich Norms on Spaces of Vector Measures, Results Math., 70
(2016), 349-371.

[7] N. Dinculeanu, Vector Measures, VEB Deutscher Verlag der
Wissenschafthen, 1966.

[8] N. Dunford, J.T. Schwartz, Linear Operators, Part I: General Theory,
Interscience Publishers, 1957.

[9] K. Falconer, Fractal Geometry. Mathematical Foundations and Applications
(third edition), Wiley, 2014.

[10] K. J. Falconer, T. C. O'Neil, Vector-valued multifractal measures,
Proc. R. Soc. Lond., Ser. A, 452 (1996), 1433-1457.

[11] P. R. Halmos, Measure Theory (eleventh printing), D. Van Nostrand
Company, 1966.

[12] J. Hutchinson, Fractals and Self-Similarity, Indiana Univ. Math. J., 30
(1981), 713-747

[13] J.L. Kelley, General topology, D. Van Nostrand Company, 1955.

[14] D. La Torre, F. Mendivil, The Monge-Kantorovich metric on multimeasures
and self-similar multimeasures, Set-Valued Var. Anal., 23 (2015), 319-331.

[15] F. Mendivil, E. R. Vrscay, Fractal vector measures and vector calculus
on planar fractal domains, Chaos Solitons Fractals, 14 (2002), 1239--1254.

[16] F. Mendivil, E. R. Vrscay, Self-affine vector measures and vector
calculus, in Fractals in Multimedia (M. F. Barnsley, D. Saupe, E. R. Vrscay,
Editors). The IMA Volumes in Mathematics and ITS Applications, vol. 132,
137-155, Springer, 2002.

[17] R. Miculescu, Generalized Iterated Function Systems with Place
Dependent Probabilities, Acta Appl. Math., 130 (2014), 135-150.

[18] A. Mihail, R. Miculescu, A Generalization of the Hutchinson Measure,
Mediterr. J. Math., 6 (2009), 203-213.

[19] N. A. Secelean, Countable Iterated Function Systems,LAP LAMBERT
Academic Publishing, 2013.

\bigskip

Ion Chi\c{t}escu

Faculty of Mathematics and Computer Science, University of Bucharest,
Academiei Str. 14, 010014, Bucharest, Romania, Email: ionchitescu@yahoo.com

\bigskip

Loredana Ioana

Faculty of Mathematics and Computer Science, University of Bucharest,
Academiei Str. 14, 010014, Bucharest, Romania, Email:
loredana.madalina.ioana@gmail.com

\bigskip

Radu Miculescu

Faculty of Mathematics and Computer Science, University of Bucharest,
Academiei Str. 14, 010014, Bucharest, Romania, Email: miculesc@yahoo.com

\bigskip

Lucian Ni\c{t}\u{a}

Technical University of Civil Engineering, Lacul Tei Blvd., 122-124, 020396,
Bucharest, Romania, Email: luci6691@yahoo.com

\end{document}